\documentclass[10pt,twoside,fleqn]{article}
\usepackage[margin=2cm]{geometry}
\usepackage{graphicx}
\usepackage{multirow}
\usepackage [displaymath, mathlines]{lineno}
\usepackage{titlesec}
\usepackage{color}
\usepackage{hyperref}
\hypersetup{
  colorlinks, linkcolor=blue, urlcolor=blue, citecolor=blue
}
\titlelabel{\thetitle.\quad}
\usepackage[font=small,format=plain,labelfont=bf,up]{caption}
\usepackage{amsmath}
\date{}
\usepackage{amsfonts}
\begin{document}
{\renewcommand{\arraystretch}{0.65}
\begin{table}[ht]
\begin{tabular}{ll}

                              
\end{tabular}
\end{table}}
\centerline{}
\centerline{}
\centerline{}


\centerline {\huge{\bf  On the Number of Cycles in a Graph}}

\centerline{}


\centerline{}

\newcommand{\mvec}[1]{\mbox{\bfseries\itshape #1}}

\centerline{\bf {Nazanin Movarraei$^{*}$}}

\centerline{}
{\small
\centerline{\emph{Department of Mathematics, University of Pune, Pune 411 007 (India) }}
\centerline{\emph{*Corresponding author E-mail:  nazanin.movarraei@gmail.com}}}

\centerline{}
\centerline{}


\newtheorem{Theorem}{Theorem}[section]

\newtheorem{Definition}[Theorem]{Definition}

\newtheorem{Proposition}[Theorem]{Proposition}

\newtheorem{Proof}[Theorem]{Proof}

\newtheorem{Example}[Theorem]{Example}

\newtheorem{Lemma}[Theorem]{Lemma}


\noindent\hrulefill

\hspace{-15 pt}{\textbf{Abstract}\\}
\centerline{}
In this paper, we obtain explicit formulae for the number of 7-cycles and the total number of paths of lengths 6 and 7 those contain a specific vertex $v_{i}$ in a simple graph G, in terms of the adjacency matrix and with the help of combinatorics.\\
\hspace{-15 pt}{{\footnotesize \emph{\textbf{Keywords}}:  \emph{ Adjacency Matrix, Cycle, Graph Theory, Path, Subgraph, Walk .}}

\noindent\hrulefill
\section{Introduction}
In a simple graph G, a walk is a sequence of vertices and edges of the form $v_{0} , e_{1} , v_{1} , ... , e_{k} , v_{k}$ such that the edge $e_{i}$ has ends $v_{i-1}$ and $v_{i}$. A walk is called closed if $v_{0} = v_{k}$. If the vertices of a walk are distinct then the walk is called a path. A cycle is a non-trivial closed walk in which all the vertices are distinct except the end vertices.\\
It is known that if a graph G has adjacency matrix A=[$a_{ij}$], then for k = 0, 1, ... , the ij-entry of A$^{k}$ is the number of $v_{i} -v_{j}$ walks of length k in G. 
It is also known that tr(A$^n$) is the sum of the diagonal entries of A$^{n}$ and $d_{i}$ is the degree of the vertex $v_{i}$.\\
In 1971, Frank Harary and Bennet Manvel [4], gave formulae for the number of cycles of lengths 3 and 4 in simple graphs as given by the following theorems:
\begin{Theorem}
 $[4]$ If G is a simple graph with adjacency matrix A, then the  number of $3$-cycles in G is $\frac{1}{6}$ tr(A$^{3}$). \\$($It is known that tr(A$^{3})=\displaystyle{\sum_{i=1}^n}a_{ii}^{(3)}=\displaystyle{\sum_{j\neq i}}a_{ij}^{(2)}a_{ij})$.
\end{Theorem}
\begin{Theorem}
$[4]$ If G is a simple graph with adjacency matrix A, then the number of $4-$cycles in G is  \\  $\frac{1}{8} [$ tr($A^{4}$)$-2$q$-2$ $\displaystyle{\sum_{j\neq i}}$$a_{ij}^{(2)}]$, where q is the number of edges in G.\\
$($It is obvious that the above formula is also equal to  $\frac{1}{8}$ $[$tr$A^{4}-$ tr$A^{2}-2$ $\displaystyle{\sum_{j\neq i}}a_{ij}^{2}]$  $  )$ 
\end{Theorem}
\begin{Theorem}
$[4]$ If G is a simple graph with n vertices and the adjacency matrix A= $[a_{ij}]$, then the number of $5-$cycles in G is $\frac{1}{10 }[$tr$(A^{5})$$+5$ tr$(A^{3})- 5$ $\displaystyle{\sum_{i=1}^n}{d_i}{a_{ii}^{(3)}}]$
\end{Theorem}
 In $2003$, V. C. Chang and H. L. Fu  [5], found a formula for the number of  $6-$cycles in a simple graph which is stated below:\\\\
\begin{Theorem}
 $[5]$ If G is a simple graph with adjacency matrix A, then the number of $6-$cycles in G is $\frac{1}{12}  [$tr$(A^{6}) -6$ tr$(A^{4}) +$5tr$(A^{3}) -4$tr$(A^{2}) -3$$\displaystyle{\sum_{i,j=1}^p}$ $a_{ij}^{(3)}$ $+12$$\displaystyle{\sum_{i,j=1}^p}$ $a_{ij}^{(2)}$ $-3$ $\displaystyle{\sum_{i=1}^p}$ $[a_{ii}^{(3)}]^2$ $+9$ $\displaystyle{\sum_{j\neq i}}$$a_{ij}^{(2)} ( a_{ij}^{(2)}-1) a_{ij}$ $-6$ $\displaystyle{\sum_{j\neq i}}$$a_{ij}^{(2)} ( a_{ij}^{(2)}-1) (a_{ij}^{(2)} -2 )$ $-2$$\displaystyle{\sum_{i=1}^p}$$a_{ii}^{(2)} ( a_{ii}^{(2)}-1)    (a_{ii}^{(2)}-2 ) ]$, where p is the number of vertices in G.
\end{Theorem}
\noindent Their proofs are based on the following fact:\\The number of n-cycles (n$=3, 4, 5, 6)$ in a graph G is equal to $\frac{1}{2n}($tr(A$^{n})- x$ ) where $x$ is the number of closed walks of length n, which are not n-cycles.\\\\
In 1997, N. Alon, R. Yuster and U. Zwick [6], gave number of 7- cyclic  graphs. They also gave some formulae for the number of cycles of lengths 5 , which contains a specific vertex $v_{i}$ in a graph G.\\
In [6, 11, 12, 13, 15, 16, 18], we have also some bounds to estimate the total time complexity for finding or counting paths and cycles in a graph.\\

\noindent In our recent works [1, 2], we obtained some formulae to find the exact number of paths of lengths 3, 4 and 5, in a simple graph G, given below:\\
\begin{Theorem}
$[1]$ Let G be a simple graph with n vertices and the adjacency matrix A= $[a_{ij}]$. The number of paths of length $3$ in G is
$\displaystyle{\sum_{j\neq i}}$$a_{ij}^{(2)}(d_{j}-a_{ij}-1)$.
\end{Theorem}
\begin{Theorem}
$[1]$ Let G be a simple graph with n vertices and the adjacency matrix A= $[a_{ij}]$. The number of paths of length $4$ in G is 
$\displaystyle{\sum_{j\neq i}}[a_{ij}^{(4)}-2a_{ij}^{(2)}(d_{j}-a_{ij})]-\displaystyle{\sum_{i=1}^n}[(2d_{i}-1)a_{ii}^{(3)}+6\left(\begin{array}{c}d_{i}\\3\end{array}\right)]$.
\end{Theorem}
\begin {Theorem}
$[1]$ Let G be a simple graph with n vertices and the adjacency matrix A= $[a_{ij}]$. The number of paths of length $3$ in G, each of which starts from a specific vertex $v_{i}$ is $\displaystyle{\sum_{j\neq i}}$$a_{ij}^{(2)}(d_{j}-a_{ij}-1)$.
\end{Theorem}
\begin{Theorem}
$[1]$ Let G be a simple graph with n vertices and the adjacency matrix A= $[a_{ij}]$. The number of paths of length $4$ in G, each of which starts from a specific vertex $v_{i}$ is $\displaystyle{\sum_{j\neq i}}[a_{ij}^{(4)}-(d_{i}+d_{j}-3a_{ij})a_{ij}^{(2)}-(a_{ii}^{(3)}+a_{jj}^{(3)}+2\left(\begin{array}{c}d_{j}-1\\2\end{array}\right))a_{ij}]$.
\end{Theorem}
\begin{Theorem}
$[2]$ Let G be a simple graph with n vertices and  the adjacency matrix A= $[a_{ij}]$. The number of paths of length $5$ in G is  $\displaystyle{\sum_{j\neq i}}a_{ij}^{(5)}-$ $2\displaystyle{\sum_{j\neq i}}a_{ij}^{(4)}+$ $2\displaystyle{\sum_{i=1}^n}a_{ii}^{(3)}(d_{i}-2)+$ $4\displaystyle{\sum_{j\neq i}}a_{ij}^{(2)}-$ $2\displaystyle{\sum_{j\neq i}}a_{ij}^{(2)}(d_{j}-a_{ij}-1)-$ $4\displaystyle{\sum_{j\neq i}}a_{ij}^{(2)}$$\left(\begin{array}{c}d_{i}-a_{ij}-1\\2\end{array}\right)$$\\+$ $6\displaystyle{\sum_{j\neq i}}a_{ij}$$\left(\begin{array}{c}a_{ij}^{(2)}\\2\end{array}\right)-$ $2\displaystyle{\sum_{j\neq i}}$$a_{ii}^{(3)} a_{ij}^{(2)}-$ $2\displaystyle{\sum_{i=1}^n}$ $a_{ii}^{(3)}$$\left(\begin{array}{c}d_{i}-2\\2\end{array}\right)$$-$ $2\displaystyle{\sum_{i=1}^n}(a_{ii}^{(4)} - a_{ii}^{(2)} -$ $2 $$\left(\begin{array}{c}d_{i}\\2\end{array}\right)$ $-\displaystyle{\sum_{j=1, j\neq i}^n} a_{ij}^{(2)})(d_{i}-2)- \displaystyle{\sum_{j\neq i}}a_{ij}-$ 3 tr $A^{4}$+ 6 tr $A^{3}$+ 3 tr $A^{2}$.
\end{Theorem}
\begin{Theorem}
$[2]$  If G is a simple graph with n vertices and the adjacency matrix A= $[a_{ij}]$, then the number of $4-$cycles each of which contains a specific vertex $v_{i}$ of G is    $\frac{1}{2}$ $[a_{ii}^{(4)} - a_{ii}^{(2)} - 2 $ $\left(\begin{array}{c}d_{i}\\2\end{array}\right)$$ $ $-\displaystyle{\sum_{j=1, j\neq i}^n} a_{ij}^{(2)}]$.
\end {Theorem}
In [6] we can also see a formula for the number of 5-cycles each of which contains a specific vertex $v_{i}$ but, their formula has some problem in coefficients. In [3] we gave the correct formula as considered bellow :
\begin{Theorem}
$[3]$ If G is a simple graph with n vertices and the adjacency matrix  A= $[a_{ij}]$, then the number of $5-$cycles each of which contains a specific vertex $v_{i}$ of G is $\frac{1}{2}$ $ [a_{ii}^{(5)} - 5 a_{ii}^{(3)} - 2 ( d_{i} - 2 ) a_{ii}^{(3)}-2$ $\displaystyle{\sum_{ j=1, j\neq i }^n}$ $a_{ij}^{(2)} a_{ij}(d_{j}-2) - 2$ $\displaystyle{\sum_{ j=1, j\neq i}^n}$ $a_{ij}(\frac{1}{2} a_{jj}^{(3)}-a_{ij}a_{ij}^{(2)})]$.\\
\end{Theorem}
In this paper we give a formula to count the exact number of cycles of length 7 and the number of cycles of lengths 6 and 7 those contain a specific vertex $v_{i}$ in a simple graph G, in terms of the adjacency matrix of G and with the helps of combinatorics.\\
\section{Number of 7- Cycles : }
In 1997, N. Alon, R. Yuster and U. Zwick [6], gave number of 7- cyclic  graphs.  The n- cyclic graph is a graph that contain a closed walk of length n and these walks are not  necessarily a cycle. In this section we obtain a formula for the number of cycles of length 7 in a simple graph G with the helps of [6]. \\\\
\noindent{\bf Method:}
To count N in the cases those considered below, we first count tr$(A^{7})$ for the graph of first configuration. This will give us the number of all closed walks of length 7 in the corresponding graph. But, some of these walks do not pass through all the edges and vertices of that configuration and to find N in each case, we have to include in any walk, all the edges and the vertices of the corresponding subgraphs at least once. So, we delete the number of closed walks of length 7 those do not pass through all the edges and vertices. To find these kind of walks we also have to count tr$(A^{7})$ for all the subgraphs of the corresponding graph those can contain a closed walk of length 7.
\begin{Theorem}:
 If G is a simple graph with n vertices and the adjacency matrix A= $[a_{ij}]$, then the number of $7-$cycles in G is $\frac{1}{14}  $ $($tr$(A^{7})- x )$, where x is equal to $\displaystyle{\sum_{i=1}^{11}}$$F_{i}$ in the cases those are considered below.\\
\end{Theorem}
\noindent {\bf Proof:} {The number of $7-$cycles of a graph G is equal to $\frac{1}{14}$$(tr(A^{7})-x)$, where $x$ is the number of closed walks of length 7 that are not 7-cycles. To find $x$, we have 11 cases as considered below; the cases are based on the configurations-(subgraphs) that generate all closed walks of length 7 that are not 7-cycles. In each case, N denote the number of closed walks of length 7 that are not 7-cycles in the corresponding subgraph, M denote the number of subgraphs of G of the same configuration and F$_{n}$, (n= 1, 2, ...) denote the total number of closed walks of length 7 that are not cycles in all possible subgraphs of G of the same configuration. However, in the cases with more that one figure ( Cases 5, 6, 8, 9, 11), N, M and F$_{n}$ are based on the first graph of the respective figures and P$_{1}$, P$_{2}$, ... denote the number of subgraphs of G which don't have the same configuration as the first graph but are counted in M. It is clear that F$_{n}$ is equal to N$\times ~($M$-$ P$_{1}-$ P$_{2}-$ ...). To find N in each case, we have to include in any walk, all the edges and the vertices of the corresponding subgraphs at least once.\\\\
\noindent {\bf Case 1:}
For the configuration of Fig $1$, N= $126$, M= $\frac{1}{6} $tr($A^{3}$) and F$_{1}$= $21$ tr($A^{3}$). (See theorem 1.1)\\
\unitlength=0.9mm
\special{em:linewidth 0.4pt}
\linethickness{0.4pt}
\begin{picture}(83.33,26)
\put(82.67,2){\line(1,2){10}}
\put(103.33,2){\line(-1,2){10}}
\put(82.67,2){\line(1,0){20.67}}
\put(92.90,22){\circle*{1.33}}
\put(82.67,2){\circle*{1.33}}
\put(103.33,2){\circle*{1.33}}
\put(93.33,-2){\makebox(0,0)[cc]{Fig 1}}
\end{picture}\\\\
\noindent {\bf Case 2:}
For the configuration of Fig $2$, N= $84$, M= $\frac{1}{2}$ $\displaystyle{\sum_{i=1}^n}$ $a_{ii}^{(3)}(d_{i}-2)$ and F$_{2}$= $42$ $\displaystyle{\sum_{i=1}^n}$ $a_{ii}^{(3)}(d_{i}-2)$.\\
\unitlength=0.9mm
\special{em:linewidth 0.4pt}
\linethickness{0.4pt}
\begin{picture}(83.33,26)
\put(82.67,2){\line(1,2){10}}
\put(103.33,2){\line(-1,2){10}}
\put(92.90,22){\line(1,0){20.67}}
\put(82.67,2){\line(1,0){20.67}}
\put(92.90,22){\circle*{1.33}}
\put(113.90,22){\circle*{1.33}}
\put(82.67,2){\circle*{1.33}}
\put(103.33,2){\circle*{1.33}}
\put(93.33,-2){\makebox(0,0)[cc]{Fig 2}}
\end{picture}\\\\
\noindent {\bf Case 3:}
For the configuration of Fig $3$, N= $28$, M= $\frac{1}{2}$ $\displaystyle{\sum_{i=1}^n}$ $a_{ii}^{(3)}$$\left(\begin{array}{c}d_{i}-2\\2\end{array}\right)$$ $ and F$_{3}$= $14$$\displaystyle{\sum_{i=1}^n}$ $a_{ii}^{(3)}$$\left(\begin{array}{c}d_{i}-2\\2\end{array}\right)$.\\
\unitlength=0.9mm
\special{em:linewidth 0.4pt}
\linethickness{0.4pt}
\begin{picture}(83.33,30)
\put(82.67,2){\line(1,2){10}}
\put(103.33,2){\line(-1,2){10}}
\put(92.90,22){\line(1,0){20.67}}
\put(92.90,22){\line(-1,0){20.67}}
\put(82.67,2){\line(1,0){20.67}}
\put(92.90,22){\circle*{1.33}}
\put(72.90,22){\circle*{1.33}}
\put(113.90,22){\circle*{1.33}}
\put(82.67,2){\circle*{1.33}}
\put(103.33,2){\circle*{1.33}}
\put(93.33,-2){\makebox(0,0)[cc]{Fig 3}}
\end{picture}\\\\\\
\noindent {\bf Case 4:}
For the configuration of Fig $4$, N= $112$, M= $\frac{1}{2}$ $\displaystyle{\sum_{j\neq i}}$$\left(\begin{array}{c}a_{ij}^{(2)}\\2\end{array}\right)$$a_{ij} $ and F$_{4}$= $56$ $\displaystyle{\sum_{j\neq i}}$$\left(\begin{array}{c}a_{ij}^{(2)}\\2\end{array}\right)$$a_{ij} $.\\
\unitlength=0.9mm
\special{em:linewidth 0.4pt}
\linethickness{0.4pt}
\begin{picture}(83.33,28)
\put(82.67,2){\line(1,2){10}}
\put(103.33,2){\line(-1,2){10}}
\put(103.33,2){\line(1,2){10}}
\put(92.90,22){\line(1,0){20.67}}
\put(82.67,2){\line(1,0){20.67}}
\put(92.90,22){\circle*{1.33}}
\put(113.90,22){\circle*{1.33}}
\put(82.67,2){\circle*{1.33}}
\put(103.33,2){\circle*{1.33}}
\put(93.33,-2){\makebox(0,0)[cc]{Fig 4}}
\end{picture}\\\\\\
\noindent {\bf Case 5:}
{For the configuration of Fig $5(a)$, N= $14$, M= $\frac{1}{2}$$\displaystyle{\sum_{j\neq i}}$$a_{ii}^{(3)} a_{ij}^{(2)} $. Let P$_{1}$ denote the number of subgraphs of G that have the same configuration as the graph of Fig $5$(b) and are counted in M. Thus  P$_{1}$= $6\times$$(\frac{1}{126}$ F$_{1})$, where $\frac{1}{126}$ F$_{1}$ is the number of subgraphs of G that have the same configuration as the graph of  Fig $5$(b) and $6$ is the number of times that this subgraph is counted in M. Let P$_{2}$ denote the number of subgraphs of G that have the same configuration as the graph of Fig 5(c) and are counted in M. Thus P$_{2}$= $2\times$$(\frac{1}{84}$ F$_{2})$, where $\frac{1}{84}$ F$_{2}$ is the number of subgraphs of G that have the same configuration as the graph of Fig $5$(c) and $2$ is the number of times that this subgraph is counted in M. Let P$_{3}$ denote the number of subgraphs of G that have the same configuration as the graph of Fig 5(d) and are counted in M. Thus P$_{3}$= $4\times$$(\frac{1}{112}$ F$_{4})$, where $\frac{1}{112}$ F$_{4}$ is the number of subgraphs of G that have the same configuration as the graph of Fig $5$(d) and $4$ is the number of times that this subgraph is counted in M. Consequently, F$_{5}$= $7$$\displaystyle{\sum_{j\neq i}}$ $a_{ii}^{(3)} a_{ij}^{(2)}- $ $\frac{2}{3}$ F$_{1}-$  $\frac{1}{3}$ F$_{2}-$ $\frac{1}{2}$ F$_{4}$ .\\
\unitlength=0.9mm
\special{em:linewidth 0.4pt}
\linethickness{0.4pt}
\begin{picture}(136,30)
\put(118.33,2){\line(-1,2){10}}
\put(107.90,22){\line(1,0){20.67}}
\put(27,22){\line(-1,-2){10}}
\put(97.90,2){\line(1,2){10}}
\put(47,22){\line(-1,0){20.67}}
\put(97.67,2){\line(1,0){20.67}}
\put(26.90,22){\circle*{1.33}}
\put(17,2){\circle*{1.33}}
\put(97.90,2){\circle*{1.33}}
\put(107.90,22){\circle*{1.33}}
\put(117.67,2){\circle*{1.33}}
\put(128.33,22){\circle*{1.33}}
\put(92.33,-8){\makebox(0,0)[cc]{Fig 5}}
\put(107.67,-2){\makebox(0,0)[cc]{$( c)$}}
\put(37.67,2){\line(1,2){10}}
\put(58.33,2){\line(-1,2){10}}
\put(37.67,2){\line(1,0){20.67}}
\put(47.90,22){\circle*{1.33}}
\put(37.67,2){\circle*{1.33}}
\put(58.33,2){\circle*{1.33}}
\put(47,-2){\makebox(0,0)[cc]{$( a )$}}
\unitlength=0.9mm
\special{em:linewidth 0.4pt}
\linethickness{0.4pt}
\put(132.67,2){\line(1,2){10}}
\put(153.33,2){\line(-1,2){10}}
\put(153.33,2){\line(1,2){10}}
\put(142.90,22){\line(1,0){20.67}}
\put(132.67,2){\line(1,0){20.67}}
\put(142.90,22){\circle*{1.33}}
\put(163.90,22){\circle*{1.33}}
\put(132.67,2){\circle*{1.33}}
\put(153.33,2){\circle*{1.33}}
\put(143,-2){\makebox(0,0)[cc]{$( d)$}}
\put(67.67,2){\line(1,2){10}}
\put(88.33,2){\line(-1,2){10}}
\put(67.67,2){\line(1,0){20.67}}
\put(77.90,22){\circle*{1.33}}
\put(67.67,2){\circle*{1.33}}
\put(88.33,2){\circle*{1.33}}
\put(77.33,-2){\makebox(0,0)[cc]{$(b)$}}
\end{picture}\\\\\\\\
\noindent {\bf Case 6:}
For the configuration of Fig $6(a)$, N= $14$, M= $\frac{1}{2}$ $\displaystyle{\sum_{ j\neq i }}$ $a_{ij}^{(2)} a_{ij}(d_{j}-2)(d_{i}-2)$. Let P$_{1}$ denote the number of subgraphs of G that have the same configuration as the graph of Fig $6( b)$ and are counted in M. Thus P$_{1}$=  $2\times$$(\frac{1}{112}$ F$_{4})$, where $\frac{1}{112}$ F$_{4}$ is the number of subgraphs of G that have the same configuration as the graph of Fig $6( b)$ and $2$  is the number of times that this subgraph is counted in M. Consequently, F$_{6}$= $7$ $\displaystyle{\sum_{ j\neq i }}$ $a_{ij}^{(2)} a_{ij}(d_{j}-2)(d_{i}-2)-$ $\frac{1}{4}$ F$_{4}$.\\
\unitlength=0.9mm
\special{em:linewidth 0.4pt}
\linethickness{0.4pt}
\begin{picture}(100,25)
\put(123.33,2){\line(1,2){10}}
\put(102.67,2){\line(1,2){10}}
\put(123.33,2){\line(-1,2){10}}
\put(112.90,22){\line(1,0){20.67}}
\put(102.67,2){\line(1,0){20.67}}
\put(112.90,22){\circle*{1.33}}
\put(102.67,2){\circle*{1.33}}
\put(123.33,2){\circle*{1.33}}
\put(133.33,22){\circle*{1.33}}
\put(98.33,-8){\makebox(0,0)[cc]{Fig 6}}
\put(112.67,-2){\makebox(0,0)[cc]{$( b )$}}

\put(72.67,2){\line(1,2){10}}
\put(93.33,2){\line(-1,2){10}}
\put(72.67,2){\line(-1,0){20.67}}
\put(82.90,22){\line(-1,0){20.67}}
\put(72.67,2){\line(1,0){20.67}}
\put(82.90,22){\circle*{1.33}}
\put(62.90,22){\circle*{1.33}}
\put(52.90,2){\circle*{1.33}}
\put(72.67,2){\circle*{1.33}}
\put(93.33,2){\circle*{1.33}}
\put(83,-2){\makebox(0,0)[cc]{$( a )$}}
\end{picture}\\\\\\\\
\noindent {\bf Case 7:}
For the configuration of Fig $7$, N= $70$, M= $\frac{1}{10 }[$ tr$(A^{5})$$+~ 5$ tr$(A^{3})- 5$ $\displaystyle{\sum_{i=1}^n}{d_i}{a_{ii}^{(3)}}]$ (See Theorem 1.3) and F$_{7}$= $7$tr($A^{5}$) + $35$ tr$(A^{3}$) $- $ $35$ $\displaystyle{\sum_{i=1}^n}{d_i}{a_{ii}^{(3)}}$.\\
\unitlength=0.9mm
\special{em:linewidth 0.4pt}
\linethickness{0.4pt}
\begin{picture}(83.33,25)
\put(90.67,2){\line(-1,2){7}}
\put(108.33,2){\line(1,2){7}}
\put(115.5,16.5){\line(-1,1){16}}
\put(83.90,16.5){\line(1,1){16}}
\put(90.67,2){\line(1,0){18}}
\put(83.90,16){\circle*{1.33}}
\put(99.7,32){\circle*{1.33}}
\put(115.5,16.5){\circle*{1.33}}
\put(90.67,2){\circle*{1.33}}
\put(108.33,2){\circle*{1.33}}
\put(97.33,-2){\makebox(0,0)[cc]{Fig 7}}
\end{picture}\\\\\\
\noindent {\bf Case 8:}
For the configuration of Fig $8( a )$, N= $42$, M= $\frac{1}{2}$ $\displaystyle{\sum_{j\neq i }}$ $a_{ij}^{(2)}(d_{j}-a_{ij}-1)(a_{ij} a_{ij}^{(2)})$ (See Theorem 1.5). Let P$_{1}$ denotes the number of subgraphs of G that have the same configuration as the graph of Fig $8( b )$ and are counted in M. Thus P$_{1}$= $4\times$$(\frac{1}{112}$ F$_{4})$, where $\frac{1}{112}$ F$_{4}$ is the number of subgraphs of G that have the same configuration as the graph of Fig $8(b)$ and $4$ is the number of times that this subgraph is counted in M. Consequently, \\F$_{8}$= $21$ $\displaystyle{\sum_{j\neq i }}$ $a_{ij}^{(2)}(d_{j}-a_{ij}-1)(a_{ij} a_{ij}^{(2)})-$ $\frac{3}{2}$ F$_{4}$.\\
\unitlength=0.9mm
\special{em:linewidth 0.4pt}
\linethickness{0.4pt}
\begin{picture}(83.33,32)
\put(75.67,2){\line(-1,2){7}}
\put(93.33,2){\line(1,2){7}}
\put(100.5,16.5){\line(-1,1){16}}
\put(68.90,16.5){\line(1,1){16}}
\put(75.67,2){\line(1,0){18}}
\put(68.90,16){\line(1,0){31}}
\put(68.90,16){\circle*{1.33}}
\put(84.7,32){\circle*{1.33}}
\put(100.5,16.5){\circle*{1.33}}
\put(75.67,2){\circle*{1.33}}
\put(93.33,2){\circle*{1.33}}
\put(97.33,-6){\makebox(0,0)[cc]{Fig 8}}
\put(112.67,-2){\makebox(0,0)[cc]{$(b)$}}
\put(84.67,-2){\makebox(0,0)[cc]{$($a$)$}}
\put(123.33,2){\line(1,2){10}}
\put(102.67,2){\line(1,2){10}}
\put(123.33,2){\line(-1,2){10}}
\put(112.90,22){\line(1,0){20.67}}
\put(102.67,2){\line(1,0){20.67}}
\put(112.90,22){\circle*{1.33}}
\put(102.67,2){\circle*{1.33}}
\put(123.33,2){\circle*{1.33}}
\put(133.33,22){\circle*{1.33}}
\end{picture}\\\\\\
\noindent {\bf Case 9:}
For the configuration of Fig $9( a )$, N= $14$, M= $\frac{1}{2}$ $[\displaystyle{\sum_{i=1}^n}$($a_{ii}^{(5)} - 5 a_{ii}^{(3)} - 2 ( d_i - 2 ) a_{ii}^{(3)})(d_{i} - 2)$$-2$$\displaystyle{\sum_{ j\neq i }}$ $a_{ij}^{(2)} a_{ij}(d_{j}-2)(d_{i}-2)$ $- 2$ $\displaystyle{\sum_{j\neq i }}$ $a_{ij}(d_{i}-2)(\frac{1}{2} a_{jj}^{(3)}-a_{ij}a_{ij}^{(2)})]$ (See Theorem 1.11). Let P$_{1}$ denote the number of subgraphs of G that have the same configuration as the graph of Fig $9( b )$ and are counted in M. Thus P$_{1}$= $2\times(\frac{1}{42}$ F$_{8})$, where $\frac{1}{42}$ F$_{8}$ is the number subgraphs of G that have the same configuration as the graph of  Fig $9(b)$ and $2$ is the number of times that this subgraph is counted in  M. Consequently, F$_{9}$= $7$ $\displaystyle{\sum_{i=1}^n}$($a_{ii}^{(5)} - 5 a_{ii}^{(3)} - 2 ( d_i - 2 ) a_{ii}^{(3)})(d_{i} - 2)-14$$\displaystyle{\sum_{ j\neq i }}$ $a_{ij}^{(2)} a_{ij}(d_{j}-2)(d_{i}-2)$ $- 14$ $\displaystyle{\sum_{j\neq i }}$ $a_{ij}(d_{i}-2)(\frac{1}{2} a_{jj}^{(3)}-a_{ij}a_{ij}^{(2)})-$ $\frac{2}{3}$ F$_{8}$.\\
\unitlength=0.9mm
\special{em:linewidth 0.4pt}
\linethickness{0.4pt}
\begin{picture}(83.33,32)
\put(75.67,2){\line(-1,2){7}}
\put(93.33,2){\line(1,2){7}}
\put(100.5,16.5){\line(-1,1){16}}
\put(68.90,16.5){\line(1,1){16}}
\put(75.67,2){\line(1,0){18}}
\put(84.7,32){\line(-1,0){22}}
\put(68.90,16){\circle*{1.33}}
\put(84.7,32){\circle*{1.33}}
\put(100.5,16.5){\circle*{1.33}}
\put(75.67,2){\circle*{1.33}}
\put(93.33,2){\circle*{1.33}}
\put(63.48,32){\circle*{1.33}}
\put(104.33,-6){\makebox(0,0)[cc]{Fig 9}}
\put(123.67,-2){\makebox(0,0)[cc]{$(b)$}}
\put(84.67,-2){\makebox(0,0)[cc]{$($a$)$}}
\put(115.67,2){\line(-1,2){7}}
\put(133.33,2){\line(1,2){7}}
\put(115.67,2){\line(1,3){10}}
\put(140.5,16.5){\line(-1,1){16}}
\put(108.90,16.5){\line(1,1){16}}
\put(115.67,2){\line(1,0){18}}
\put(108.90,16){\circle*{1.33}}
\put(125.2,31.8){\circle*{1.33}}
\put(140.5,16.5){\circle*{1.33}}
\put(115.67,2){\circle*{1.33}}
\put(133.33,2){\circle*{1.33}}
\end{picture}\\\\\\\\
\noindent {\bf Case 10:}
For the configuration of Fig $10$, N= $84$, M= $\frac{1}{2}$ $\displaystyle{\sum_{j\neq i}} $$\left(\begin{array}{c}a_{ij}^{(2)}\\3\end{array}\right)$$ a_{ij}$ and F$_{10}$=$~42$$\displaystyle{\sum_{j\neq i}} $$\left(\begin{array}{c}a_{ij}^{(2)}\\3\end{array}\right)$$ a_{ij}$.\\
\unitlength=0.9mm
\special{em:linewidth 0.4pt}
\linethickness{0.4pt}
\begin{picture}(83.33,35)
\put(92.67,2){\line(1,0){20}}
\put(112.67,2){\line(-1,1){20}}
\put(92.67,2){\line(0,1){20}}
\put(112.67,2){\line(0,1){20}}
\put(112.67,2){\line(1,3){10}}
\put(92.67,22){\line(1,0){20}}
\put(92.67,22){\line(3,1){30}}
\put(92.67,2){\circle*{1.33}}
\put(112.67,2){\circle*{1.33}}
\put(112.67,22){\circle*{1.33}}
\put(92.67,22){\circle*{1.33}}
\put(122.3,31.6){\circle*{1.33}}
\put(101.90,-4){\makebox(0,0)[cc]{Fig 10}}
\end{picture}\\\\\\
\noindent {\bf Case 11:}
For the configuration of Fig $11( a )$, N= $28$, M= $\frac {1}{2}$ $\displaystyle{\sum_{j\neq i}} \left(\begin{array}{c}a_{ij}^{(2)}\\2\end{array}\right)a_{ii}^{(3)}$. Let P$_{1}$ denote the number of subgraphs of G that have the same configuration as the graph of Fig $11 ( b )$ and are counted in M. Thus  P$_{1}$= $2$ $\times$$(\frac{1}{42}$ F$_{8})$, where $\frac{1}{42}$ F$_{8}$ is the number of subgraphs of G that have the same configuration as the graph of Fig $11(b)$ and $2$ is the number of times that this subgraph is counted in M. Let P$_{2}$ denote the number of subgraphs of G that have the same configuration as the graph of Fig $11 ( c )$ and are counted in M. Thus P$_{2} =$ $6\times$ $(\frac{1}{112}$ F$_{4})$, where $\frac{1}{112}$ F$_{4}$ is the number of subgraphs of G that have the same configuration as the graph of Fig $11(c)$ and 6 is the number of times that this subgraph is counted in M. Let P$_{3}$ denote the number of subgraphs of G that have the same configuration as the graph of Fig $11 ( d )$ and are counted in M. Thus P$_{3}$= $6\times$ $(\frac{1}{84}$ F$_{10})$, where $\frac{1}{84}$ F$_{10}$ is the number of subgraphs of G that have the same configuration as the graph of Fig $11(d)$ and $6$ is the number of times that this subgraph is counted in M. Consequently, F$_{11}$= $14$ $\displaystyle{\sum_{j\neq i}} \left(\begin{array}{c}a_{ij}^{(2)}\\2\end{array}\right)a_{ii}^{(3)}-$ $\frac{4}{3}$ F$_{8}-$ $\frac{3}{2}$ F$_{4}-$ $2$ F$_{10}$ .\\
\unitlength=0.9mm
\special{em:linewidth 0.4pt}
\linethickness{0.4pt}
\begin{picture}(83.33,35)
\put(72.67,2){\line(1,0){20}}
\put(72.67,2){\line(0,1){20}}
\put(92.67,2){\line(0,1){20}}
\put(92.67,22){\line(-1,1){10}}
\put(72.67,22){\line(1,1){10}}
\put(72.67,22){\line(1,0){20}}
\put(72.67,2){\circle*{1.33}}
\put(82.67,32){\circle*{1.33}}
\put(92.67,2){\circle*{1.33}}
\put(92.67,22){\circle*{1.33}}
\put(72.67,22){\circle*{1.33}}
\put(81.90,-4){\makebox(0,0)[cc]{$(b)$}}
\put(32.67,2){\line(1,0){20}}
\put(32.67,2){\line(0,1){20}}
\put(52.67,2){\line(0,1){20}}
\put(32.67,22){\line(1,0){20}}
\put(32.67,22){\line(0,1){10}}
\put(32.67,22){\line(-1,0){10}}
\put(32.67,32){\line(-1,-1){10}}
\put(32.67,2){\circle*{1.33}}
\put(52.67,2){\circle*{1.33}}
\put(22.67,22){\circle*{1.33}}
\put(32.67,32){\circle*{1.33}}
\put(52.67,22){\circle*{1.33}}
\put(32.67,22){\circle*{1.33}}
\put(42.90,-4){\makebox(0,0)[cc]{$($a$)$}}
\put(112.67,2){\line(1,0){20}}
\put(132.67,2){\line(-1,1){20}}
\put(112.67,2){\line(0,1){20}}
\put(132.67,2){\line(0,1){20}}
\put(112.67,22){\line(1,0){20}}
\put(112.67,2){\circle*{1.33}}
\put(132.67,2){\circle*{1.33}}
\put(132.67,22){\circle*{1.33}}
\put(112.67,22){\circle*{1.33}}
\put(121.90,-4){\makebox(0,0)[cc]{$(c)$}}
\put(152.67,2){\line(1,0){20}}
\put(172.67,2){\line(-1,1){20}}
\put(152.67,2){\line(0,1){20}}
\put(172.67,2){\line(0,1){20}}
\put(172.67,2){\line(1,3){10}}
\put(152.67,22){\line(1,0){20}}
\put(152.67,22){\line(3,1){30}}
\put(152.67,2){\circle*{1.33}}
\put(172.67,2){\circle*{1.33}}
\put(172.67,22){\circle*{1.33}}
\put(152.67,22){\circle*{1.33}}
\put(182.3,31.6){\circle*{1.33}}
\put(161.90,-4){\makebox(0,0)[cc]{$(d)$}}
\put(101.90,-8){\makebox(0,0)[cc]{Fig 11}}
\end{picture}\\\\\\\\
Now we add the values of F$_{n}$ arising from the above cases and determine $x$. By putting the value of $x$ in  $\frac{1}{14}$$($tr$(A^{7})-x)$ we get the desired formula.\hfill$\Box$
\begin{Example}
In $K_7$, $F_{1}$= $4410$, $F_{2}$= $35280$, $F_{3}$= $17640$, $F_{4}$= $23520$, $F_{5}$= $17640$, $F_{6}$= $17640$, $F_{7}$= $17640$, \\$F_{8}$= $52920$, $F_{9}$= $35280$, $F_{10}$= $17640$, $F_{11}$= $35280$ and $tr A^{7}$= $279930$. So x= $274890$ and by Theorem $2.1$, the number of cycles of length $7$ in  K$_7$ is $\frac{1}{14}(279930-274890)=360$.
\end{Example}
\section{Number of Cycles Passing the Vertex V$_{i}$ :}
{In this section we give formulae to count the number of cycles of lengths 6 and 7, each of which contain a specific vertex v$_{i}$ of the graph G.
\begin{Theorem}
If G is a simple graph with n vertices and the adjacency matrix  A= $[a_{ij}]$, then the number of $6-$cycles each of which contains a specific vertex $v_{i}$ of G is $\frac{1}{2}$ $ (a_{ii}^{(6)}-$ x$)$, where x is equal to $\displaystyle{\sum_{i=1}^{17}}$$F_{i}$ in the cases those are considered below.
\end{Theorem}
\noindent {\bf Proof:} {The number of $6-$cycles each of which contain a specific vertex $v_{i}$ of the graph G is equal to $\frac{1}{2}$ $(a_{ii}^{(6)}-x)$, where $x$ is the number of closed walks of length $6$ form the vertex $v_{i}$ to $v_{i}$ that are not $6-$cycles. To find $x$, we have 17 cases as considered below; the cases are based on the configurations-(subgraphs) that generate $v_{i}-v_{i}$ walks of length 6 that are not cycles. In each case, N denote the number of walks of length $6$ from $v_{i}$ to $v_{i}$ that are not cycles in the corresponding subgraph, M denote the number of subgraphs of G of the same configuration and  F$_{n}$, (n= 1, 2, ...)  denote the total number of $v_{i}-v_{i}$ walks of length $6$ that are not cycles in all possible subgraphs of G of the same configuration. However, in the cases with more than one figure (Cases 11, 12, 13, 14, 15, 16, 17), N, M and F$_{n}$ are based on the first graph of the respective figures and P$_{1}$, P$_{2}$,... denote the number of subgraphs of G which don't have the same configuration as the first graph but are counted in M. It is clear that  F$_{n}$ is equal to N $\times$ (M$-$ P$_{1}-$P$_{2}-$...). To find N in each case, we have to include in any walk, all the edges and the vertices of the corresponding subgraphs at least once.\\\\
\noindent {\bf Case 1:}
{For the configuration of Fig 12, N= 1, M= $a_{ii}^{(2)}$ and F$_{1}$= $a_{ii}^{(2)}$.\\
\unitlength=0.9mm
\special{em:linewidth 0.4pt}
\linethickness{0.4pt}
\begin{picture}(90,14)
\put(103.33,7){\line(-1,0){20.67}}
\put(82.67,7){\circle*{1.33}}
\put(103.33,7){\circle*{1.33}}
\put(93.33,3){\makebox(0,0)[cc]{Fig 12}}
\put(82.67,9){\makebox(0,0)[cc]{$v_{i}$}}
\end{picture}\\\\
\noindent {\bf Case 2:}
{For the configuration of Fig 13, N= 3, M=$\displaystyle{\sum_{j\neq i}}$ $a_{ij}^{(2)}$ and F$_{2}$= 3$\displaystyle{\sum_{j\neq i}}$ $a_{ij}^{(2)}$.\\\\
\unitlength=0.9mm
\special{em:linewidth 0.4pt}
\linethickness{0.4pt}
\begin{picture}(90,14)
\put(93.33,7){\line(-1,0){20.67}}
\put(93.33,7){\line(1,0){20.67}}
\put(72.67,7){\circle*{1.33}}
\put(93.33,7){\circle*{1.33}}
\put(114,7){\circle*{1.33}}
\put(93.33,3){\makebox(0,0)[cc]{Fig 13}}
\put(72.67,9){\makebox(0,0)[cc]{$v_{i}$}}
\end{picture}\\\\
\noindent {\bf Case 3:}
{For the configuration of Fig 14,  N= 6, M=$~$$\left(\begin{array}{c}d_{i}\\2\end{array}\right)$$ $ and F$_{3}$=$~6 $$\left(\begin{array}{c}d_{i}\\2\end{array}\right)$$ $.\\
\unitlength=0.9mm
\special{em:linewidth 0.4pt}
\linethickness{0.4pt}
\begin{picture}(90,14)
\put(93.33,7){\line(-1,0){20.67}}
\put(93.33,7){\line(1,0){20.67}}
\put(72.67,7){\circle*{1.33}}
\put(93.33,7){\circle*{1.33}}
\put(114,7){\circle*{1.33}}
\put(93.33,3){\makebox(0,0)[cc]{Fig 14}}
\put(93.33,9){\makebox(0,0)[cc]{$v_{i}$}}
\end{picture}\\\\
\noindent {\bf Case 4:}
{For the configuration of Fig 15, N= 1, M=$\displaystyle{\sum_{j\neq i}}a_{ij}^{(2)}(d_{j}-a_{ij}-1)$ and F$_{4}$=$\displaystyle{\sum_{j\neq i}}a_{ij}^{(2)}(d_{j}-a_{ij}-1)$.\\ (See Theorem 1.7)\\
\unitlength=0.9mm
\special{em:linewidth 0.4pt}
\linethickness{0.4pt}
\begin{picture}(90,15)
\put(103.33,10){\line(-1,0){20.67}}
\put(83,10){\line(-1,0){20.67}}
\put(103.33,10){\line(1,0){20.67}}
\put(82.67,10){\circle*{1.33}}
\put(103.33,10){\circle*{1.33}}
\put(124,10){\circle*{1.33}}
\put(62,10){\circle*{1.33}}
\put(93.33,5){\makebox(0,0)[cc]{Fig 15}}
\put(62,13){\makebox(0,0)[cc]{$v_{i}$}}
\end{picture}\\
\noindent {\bf Case 5:}
{For the configuration of Fig 16, N= 2, M=$\displaystyle{\sum_{j\neq i}}a_{ij}^{(2)}(d_{i}-a_{ij}-1)$ and F$_{5}$= 2$\displaystyle{\sum_{j\neq i}}a_{ij}^{(2)}(d_{i}-a_{ij}-1)$.\\
\unitlength=0.9mm
\special{em:linewidth 0.4pt}
\linethickness{0.4pt}
\begin{picture}(90,19)
\put(103.33,10){\line(-1,0){20.67}}
\put(83,10){\line(-1,0){20.67}}
\put(103.33,10){\line(1,0){20.67}}
\put(82.67,10){\circle*{1.33}}
\put(103.33,10){\circle*{1.33}}
\put(124,10){\circle*{1.33}}
\put(62,10){\circle*{1.33}}
\put(93.33,5){\makebox(0,0)[cc]{Fig 16}}
\put(82.67,13){\makebox(0,0)[cc]{$v_{i}$}}
\end{picture}\\
\noindent {\bf Case 6:}
{For the configuration of Fig 17, N= 2, M=$\displaystyle{\sum_{j\neq i}}a_{ij}$$\left(\begin{array}{c}d_{j}-1\\2\end{array}\right)$$ $ and F$_{6}$= 2$\displaystyle{\sum_{j\neq i}}a_{ij}$$\left(\begin{array}{c}d_{j}-1\\2\end{array}\right)$$ $.\\
\unitlength=0.9mm
\special{em:linewidth 0.4pt}
\linethickness{0.4pt}
\begin{picture}(83.33,34)
\put(80.67,4){\line(3,2){12}}
\put(92.5,12){\line(0,1){15}}
\put(104.33,4){\line(-3,2){12}}
\put(92.7,12){\circle*{1.33}}
\put(80.67,4){\circle*{1.33}}
\put(104.33,4){\circle*{1.33}}
\put(92.5,27){\circle*{1.33}}
\put(93,2){\makebox(0,0)[cc]{Fig 17}}
\put(93,30){\makebox(0,0)[cc]{$v_{i}$}}
\end{picture}\\\\
\noindent {\bf Case 7:} {For the configuration of Fig 18, N= 6, M= $\left(\begin{array}{c}d_{i}\\3\end{array}\right)$ and F$_{7}$=  6$\left(\begin{array}{c}d_{i}\\3\end{array}\right)$.\\
\unitlength=0.9mm
\special{em:linewidth 0.4pt}
\linethickness{0.4pt}
\begin{picture}(83.33,30)
\put(80.67,4){\line(3,2){12}}
\put(92.5,12){\line(0,1){15}}
\put(104.33,4){\line(-3,2){12}}
\put(92.7,12){\circle*{1.33}}
\put(80.67,4){\circle*{1.33}}
\put(104.33,4){\circle*{1.33}}
\put(92.5,27){\circle*{1.33}}
\put(93,2){\makebox(0,0)[cc]{Fig 18}}
\put(89,13){\makebox(0,0)[cc]{$v_{i}$}}
\end{picture}\\\\

\noindent {\bf Case 8:}
{For the configuration of  Fig 19, N= 8, M= $\frac{1}{2}\displaystyle{\sum_{j\neq i}}a_{ij}^{(2)}a_{ij}$ and F$_{8}$= $4\displaystyle{\sum_{j\neq i}}a_{ij}^{(2)}a_{ij}$.\\
\unitlength=0.9mm
\special{em:linewidth 0.4pt}
\linethickness{0.4pt}
\begin{picture}(83.33,30)
\put(82.67,2){\line(1,2){10}}
\put(103.33,2){\line(-1,2){10}}
\put(82.67,2){\line(1,0){20.67}}
\put(93,22){\circle*{1.33}}
\put(82.67,2){\circle*{1.33}}
\put(103.33,2){\circle*{1.33}}
\put(93.5,25){\makebox(0,0)[cc]{$v_{i}$}}
\put(93.33,-2){\makebox(0,0)[cc]{Fig 19}}
\end{picture}\\\\
\noindent {\bf Case 9:}
For the configuration of Fig 20, N= 12, M=$\displaystyle{\sum_{j\neq i}}$$\left(\begin{array}{c}a_{ij}^{(2)}\\2\end{array}\right)$ and F$_{9}$= 12$\displaystyle{\sum_{j\neq i}}$$\left(\begin{array}{c}a_{ij}^{(2)}\\2\end{array}\right)$.\\
\unitlength=0.9mm
\special{em:linewidth 0.4pt}
\linethickness{0.4pt}
\begin{picture}(83.33,32)
\put(82.67,6){\line(1,0){20.67}}
\put(82.67,6){\line(0,1){20.67}}
\put(103,6){\line(0,1){20.67}}
\put(82.67,26){\line(1,0){20.67}}
\put(82.67,6){\circle*{1.33}}
\put(103,6){\circle*{1.33}}
\put(103,26){\circle*{1.33}}
\put(82.67,26){\circle*{1.33}}
\put(93.33,2){\makebox(0,0)[cc]{Fig 20}}
\put(102.67,29){\makebox(0,0)[cc]{$v_{i}$}}
\end{picture}\\\\
\noindent {\bf Case 10:}
For the configuration of Fig 21, N= 12, M=$\displaystyle{\sum_{j\neq i}}$$\left(\begin{array}{c}a_{ij}^{(2)}\\2\end{array}\right)$$ a_{ij}$ and F$_{10}$= 12$\displaystyle{\sum_{j\neq i}}$$\left(\begin{array}{c}a_{ij}^{(2)}\\2\end{array}\right)$$ a_{ij}$.\\
\unitlength=0.9mm
\special{em:linewidth 0.4pt}
\linethickness{0.4pt}
\begin{picture}(83,30)
\put(82.67,2){\line(1,2){10}}
\put(103.33,2){\line(-1,2){10}}
\put(103.33,2){\line(1,2){10}}
\put(92.90,22){\line(1,0){20.67}}
\put(82.67,2){\line(1,0){20.67}}
\put(92.90,22){\circle*{1.33}}
\put(113.90,22){\circle*{1.33}}
\put(82.67,2){\circle*{1.33}}
\put(103.33,2){\circle*{1.33}}
\put(93.33,-2){\makebox(0,0)[cc]{Fig 21}}
\put(92.90,25){\makebox(0,0)[cc]{$v_{i}$}}
\end{picture}\\\\\\
\noindent {\bf Case 11:}
For the configuration of Fig 22(a), N= 6, M= $\frac{1}{2}\displaystyle{\sum_{k\neq j,\ j,k\neq i}}a_{jk}^{(2)}a_{ij}a_{ik}a_{jk}$. Let $P_{1}$  denote the number of all subgraphs of G that have the same configuration as the graph of Fig 22(b) and are counted in M. Thus P$_{1}$= $1\times(\frac{1}{8}$ F$_{8})$, where $\frac{1}{8}$ F$_{8}$ is the number of subgraphs of G that have the same configuration as the graph of Fig 22(b) and this subgraph is counted only once in M. Consequently, F$_{11}$$= 3 \displaystyle{\sum_{k\neq j,\ j,k\neq i}}a_{jk}^{(2)}a_{ij}a_{ik}a_{jk}-$ $ \frac{3}{4}$ F$_{8}$ .\\

\unitlength=0.9mm
\special{em:linewidth 0.4pt}
\linethickness{0.4pt}
\begin{picture}(83,27)
\put(62.67,2){\line(1,2){10}}
\put(83.33,2){\line(-1,2){10}}
\put(83.33,2){\line(1,2){10}}
\put(72.90,22){\line(1,0){20.67}}
\put(62.67,2){\line(1,0){20.67}}
\put(72.90,22){\circle*{1.33}}
\put(93.90,22){\circle*{1.33}}
\put(62.67,2){\circle*{1.33}}
\put(83.33,2){\circle*{1.33}}
\put(93.33,-4){\makebox(0,0)[cc]{Fig 22}}
\put(92.90,25){\makebox(0,0)[cc]{$v_{i}$}}

\put(102.67,2){\line(1,2){10}}
\put(123.33,2){\line(-1,2){10}}
\put(102.67,2){\line(1,0){20.67}}
\put(113,22){\circle*{1.33}}
\put(102.67,2){\circle*{1.33}}
\put(123.33,2){\circle*{1.33}}
\put(113.5,25){\makebox(0,0)[cc]{$v_{i}$}}
\put(73.33,-2){\makebox(0,0)[cc]{($a$)}}
\put(113.33,-2){\makebox(0,0)[cc]{(b)}}
\end{picture}\\\\\\
\noindent {\bf Case 12:}
For the configuration of Fig 23(a), N= 2, M=$\displaystyle{\sum_{k\neq j,\ j,k\neq i}^n}\left(\begin{array}{c}a_{jk}^{(2)}\\2\end{array}\right)a_{ij}$ . Let P$_{1}$  denote the number of all subgraphs of G that have the same configuration as the graph of Fig 23(b) and are counted in M. Thus P$_{1} =2\times(\frac{1}{12}$ F$_{9})$, where $\frac{1}{12}$ F$_{9}$ is the number of subgraphs of G that have the same configuration as the graph of Fig 23(b) and 2 is the number of times that this subgraph is counted in M. Consequently, F$_{12}$=  2$\displaystyle{\sum_{k\neq j,\ j,k\neq i}^n}\left(\begin{array}{c}a_{jk}^{(2)}\\2\end{array}\right)a_{ij}- $ $ \frac{1}{3}$ F$_{9}$. \\
\unitlength=0.9mm
\special{em:linewidth 0.4pt}
\linethickness{0.4pt}
\begin{picture}(83.33,42)
\put(62.67,6){\line(1,0){20.67}}
\put(62.67,6){\line(0,1){20.67}}
\put(83,6){\line(0,1){20.67}}
\put(62.67,26){\line(1,0){20.67}}
\put(83,26){\line(1,1){14}}
\put(62.67,6){\circle*{1.33}}
\put(83,6){\circle*{1.33}}
\put(83,26){\circle*{1.33}}
\put(62.67,26){\circle*{1.33}}
\put(97,40){\circle*{1.33}}
\put(73.33,0){\makebox(0,0)[cc]{($a$)}}
\put(93.33,-2){\makebox(0,0)[cc]{Fig 23}}
\put(99,43){\makebox(0,0)[cc]{$v_{i}$}}
\put(102.67,6){\line(1,0){20.67}}
\put(102.67,6){\line(0,1){20.67}}
\put(123,6){\line(0,1){20.67}}
\put(102.67,26){\line(1,0){20.67}}
\put(102.67,6){\circle*{1.33}}
\put(123,6){\circle*{1.33}}
\put(123,26){\circle*{1.33}}
\put(102.67,26){\circle*{1.33}}
\put(113.33,0){\makebox(0,0)[cc]{(b)}}
\put(123,29){\makebox(0,0)[cc]{$v_{i}$}}
\end{picture}\\\\
\noindent {\bf Case 13:}
For the configuration of Fig 24(a), N= 4, M= $\displaystyle{\sum_{j\neq i}}\left(\begin{array}{c}a_{ij}^{(2)}\\2\end{array}\right)(d_{i}-2)$. Let P$_{1}$  denote the number of all subgraphs of G that have the same configuration as the graph of Fig 24(b) and are counted in M. Thus P$_{1} = 1\times(\frac{1}{12}$ F$_{10})$, where $\frac{1}{12}$ F$_{10}$ is the number of subgraphs of G that have the same configuration as the graph of Fig 24(b) and this subgraph is counted only once in M. Consequently,  F$_{13}$= $4\displaystyle{\sum_{j\neq i}}\left(\begin{array}{c}a_{ij}^{(2)}\\2\end{array}\right)(d_{i}-2)-$ $\frac{1}{3}$ F$_{10}$.\\ 
\unitlength=0.9mm
\special{em:linewidth 0.4pt}
\linethickness{0.4pt}
\begin{picture}(83.33,40)
\put(62.67,6){\line(1,0){20.67}}
\put(62.67,6){\line(0,1){20.67}}
\put(83,6){\line(0,1){20.67}}
\put(62.67,26){\line(1,0){20.67}}
\put(83,26){\line(1,1){14}}
\put(62.67,6){\circle*{1.33}}
\put(83,6){\circle*{1.33}}
\put(83,26){\circle*{1.33}}
\put(62.67,26){\circle*{1.33}}
\put(97,40){\circle*{1.33}}
\put(73.33,0){\makebox(0,0)[cc]{($a$)}}
\put(113.33,0){\makebox(0,0)[cc]{(b)}}
\put(93.33,-4){\makebox(0,0)[cc]{Fig 24}}
\put(123,29){\makebox(0,0)[cc]{$v_{i}$}}
\put(102.67,6){\line(1,0){20.67}}
\put(102.67,6){\line(0,1){20.67}}
\put(123,6){\line(0,1){20.67}}
\put(102.67,26){\line(1,0){20.67}}
\put(123,26){\line(-1,-1){20}}
\put(102.67,6){\circle*{1.33}}
\put(123,6){\circle*{1.33}}
\put(123,26){\circle*{1.33}}
\put(102.67,26){\circle*{1.33}}
\put(82,29){\makebox(0,0)[cc]{$v_{i}$}}
\end{picture}\\\\\\
\noindent {\bf Case 14:}
For the configuration of Fig 25(a), N= 2, M= $\displaystyle{\sum_{j\neq i}}\left(\begin{array}{c}a_{ij}^{(2)}\\2\end{array}\right)(d_{j}-2)$. Let P$_{1}$  denote the number of all subgraphs of G that have the same configuration as the graph of Fig 25(b) and are counted in M. Thus P$_{1} = 1\times(\frac{1}{12}$ F$_{10})$, where $\frac{1}{12}$ F$_{10}$ is the number of subgraphs of G that have the same configuration as the graph of Fig 25(b) and this subgraph is counted only once in M. Consequently,  F$_{14}$= $2\displaystyle{\sum_{j\neq i}}\left(\begin{array}{c}a_{ij}^{(2)}\\2\end{array}\right)(d_{j}-2)-$ $\frac{1}{6}$ F$_{10}$.\\
\unitlength=0.9mm
\special{em:linewidth 0.4pt}
\linethickness{0.4pt}
\begin{picture}(83.33,40)
\put(62.67,6){\line(1,0){20.67}}
\put(62.67,6){\line(0,1){20.67}}
\put(83,6){\line(0,1){20.67}}
\put(62.67,26){\line(1,0){20.67}}
\put(83,26){\line(1,1){14}}
\put(62.67,6){\circle*{1.33}}
\put(83,6){\circle*{1.33}}
\put(83,26){\circle*{1.33}}
\put(62.67,26){\circle*{1.33}}
\put(97,40){\circle*{1.33}}
\put(73.33,0){\makebox(0,0)[cc]{($a$)}}
\put(113.33,0){\makebox(0,0)[cc]{(b)}}
\put(93.33,-4){\makebox(0,0)[cc]{Fig 25}}
\put(123,29){\makebox(0,0)[cc]{$v_{i}$}}
\put(102.67,6){\line(1,0){20.67}}
\put(102.67,6){\line(0,1){20.67}}
\put(123,6){\line(0,1){20.67}}
\put(102.67,26){\line(1,0){20.67}}
\put(123,26){\line(-1,-1){20}}
\put(102.67,6){\circle*{1.33}}
\put(123,6){\circle*{1.33}}
\put(123,26){\circle*{1.33}}
\put(102.67,26){\circle*{1.33}}
\put(60,6){\makebox(0,0)[cc]{$v_{i}$}}
\end{picture}\\\\\\
\noindent {\bf Case 15:}
For the configuration of Fig 26(a), N= 2, M=$\displaystyle{\sum_{j\neq i}}a_{ij}^{(2)}(d_{j}-a_{ij}-1)a_{ij}(d_{j}-2)$. Let P$_{1}$  denote the number of all subgraphs of G that have the same configuration as the graph of Fig 26(b) and are counted in M. Thus P$_{1}= 2\times(\frac{1}{6}$ F$_{11})$, where $\frac{1}{6}$ F$_{11}$ is the number of subgraphs of G that have the same configuration as the graph of Fig 26(b) and 2 is the number of times that this subgraph is counted in M. Consequently, \\F$_{15}$= 2$\displaystyle{\sum_{j\neq i}}a_{ij}^{(2)}(d_{j}-a_{ij}-1)a_{ij}(d_{j}-2)- $ $\frac{2}{3}$ F$_{11}$.\\
\unitlength=0.9mm
\special{em:linewidth 0.4pt}
\linethickness{0.4pt}
\begin{picture}(83.33,37)
\put(62.67,6){\line(1,0){20.67}}
\put(62.67,6){\line(0,1){20.67}}
\put(83,6){\line(0,1){20.67}}
\put(62.67,26){\line(1,0){20.67}}
\put(83,26){\line(1,1){14}}
\put(62.67,6){\circle*{1.33}}
\put(83,6){\circle*{1.33}}
\put(83,26){\circle*{1.33}}
\put(62.67,26){\circle*{1.33}}
\put(97,40){\circle*{1.33}}
\put(73.33,0){\makebox(0,0)[cc]{($a$)}}
\put(113.33,0){\makebox(0,0)[cc]{(b)}}
\put(93.33,-4){\makebox(0,0)[cc]{Fig 26}}
\put(103,29){\makebox(0,0)[cc]{$v_{i}$}}
\put(102.67,6){\line(1,0){20.67}}
\put(102.67,6){\line(0,1){20.67}}
\put(123,6){\line(0,1){20.67}}
\put(102.67,26){\line(1,0){20.67}}
\put(123,26){\line(-1,-1){20}}
\put(102.67,6){\circle*{1.33}}
\put(123,6){\circle*{1.33}}
\put(123,26){\circle*{1.33}}
\put(102.67,26){\circle*{1.33}}
\put(63,29){\makebox(0,0)[cc]{$v_{i}$}}
\end{picture}\\\\\\
\noindent {\bf Case 16:}
For the configuration of Fig 27(a), N= 4, M= $\frac{1}{2}\displaystyle{\sum_{j\neq i}}a_{jj}^{(3)}a_{ij}^{(2)}a_{ij}$. Let P$_{1}$  denote the number of all subgraphs of G that have the same configuration as the graph of Fig 27(b) and are counted in M. Thus P$_{1}$= $2\times(\frac{1}{8}$ F$_{8})$, where $\frac{1}{8}$ F$_{8}$ is the number of subgraphs of G that have the same configuration as the graph of Fig 27(b) and 2 is the number of times that this subgraph is counted in M. Let P$_{2}$  denote the number of all subgraphs of G that have the same configuration as the graph of Fig 27(c) and are counted in M. Thus P$_{2}= 2\times(\frac{1}{6}$ F$_{11})$, where $\frac{1}{6}$ F$_{11}$ is the number of subgraphs of G that have the same configuration as the graph of Fig 27(c) and 2 is the number of times that this subgraph is counted in M. 
 Let P$_{3}$  denote the number of all subgraphs of G that have the same configuration as the graph of Fig 27(d) and are counted in M. Thus P$_{3} = 2\times(\frac{1}{12}$ F$_{10})$, where $\frac{1}{12}$ F$_{10}$ is the number of subgraphs of G that have the same configuration as the graph of Fig 27(d) and 2 is the number of times that this subgraph is counted in M. Consequently, F$_{16}$= $2\displaystyle{\sum_{j\neq i}}a_{jj}^{(3)}a_{ij}^{(2)}a_{ij}-   $ F$_{8}-$ $\frac{4}{3}$ F$_{11}- $ $\frac{2}{3}$ F$_{10}$. \\
\unitlength=0.9mm
\special{em:linewidth 0.4pt}
\linethickness{0.4pt}
\begin{picture}(83.33,36)
\put(22.67,6){\line(1,0){20.67}}
\put(43,6){\line(-1,1){20.67}}
\put(22.67,26){\line(1,0){20.67}}
\put(23,6){\line(1,1){20}}
\put(22.67,6){\circle*{1.33}}
\put(43,6){\circle*{1.33}}
\put(43,26){\circle*{1.33}}
\put(22.67,26){\circle*{1.33}}
\put(33.33,0){\makebox(0,0)[cc]{($a$)}}
\put(73.33,0){\makebox(0,0)[cc]{(b)}}
\put(113.33,0){\makebox(0,0)[cc]{(c)}}
\put(153.33,0){\makebox(0,0)[cc]{(d)}}
\put(93.33,-4){\makebox(0,0)[cc]{Fig 27}}
\put(63,29){\makebox(0,0)[cc]{$v_{i}$}}
\put(62.67,6){\line(0,1){20.67}}
\put(62.67,26){\line(1,0){20.67}}
\put(83,26){\line(-1,-1){20}}
\put(62.67,6){\circle*{1.33}}
\put(83,26){\circle*{1.33}}
\put(62.67,26){\circle*{1.33}}
\put(23,29){\makebox(0,0)[cc]{$v_{i}$}}

\put(103,29){\makebox(0,0)[cc]{$v_{i}$}}
\put(102.67,6){\line(1,0){20.67}}
\put(102.67,6){\line(0,1){20.67}}
\put(123,6){\line(0,1){20.67}}
\put(102.67,26){\line(1,0){20.67}}
\put(123,26){\line(-1,-1){20}}
\put(102.67,6){\circle*{1.33}}
\put(123,6){\circle*{1.33}}
\put(123,26){\circle*{1.33}}
\put(102.67,26){\circle*{1.33}}

\put(143,29){\makebox(0,0)[cc]{$v_{i}$}}
\put(142.67,6){\line(1,0){20.67}}
\put(142.67,6){\line(0,1){20.67}}
\put(163,6){\line(0,1){20.67}}
\put(142.67,26){\line(1,0){20.67}}
\put(143,26){\line(1,-1){20}}
\put(142.67,6){\circle*{1.33}}
\put(163,6){\circle*{1.33}}
\put(163,26){\circle*{1.33}}
\put(142.67,26){\circle*{1.33}}
\end{picture}\\\\\\
\noindent {\bf Case 17:}
For the configuration of Fig 28(a),  N= 8, M= $\left(\begin{array}{c}\frac{1}{2}a_{ii}^{(3)}\\2\end{array}\right)$. 
Let P$_{1}$  denote the number of all subgraphs of G that have the same configuration as the graph of Fig 28(b) and are counted in M. Thus P$_{1} = 1\times(\frac{1}{12}$ F$_{10})$, where $\frac{1}{12}$ F$_{10}$ is the number of subgraphs of G that have the same configuration as the graph of Fig 28(b) and this subgraph is counted only once in M. Consequently,  F$_{17}$= $8 \left(\begin{array}{c}\frac{1}{2}a_{ii}^{(3)}\\2\end{array}\right)- $ $\frac{2}{3}$ F$_{10}$.\\

\unitlength=0.9mm
\special{em:linewidth 0.4pt}
\linethickness{0.4pt}
\begin{picture}(83.33,30)
\put(62.67,6){\line(1,0){20.67}}
\put(83,6){\line(-1,1){20.67}}
\put(62.67,26){\line(1,0){20.67}}
\put(63,6){\line(1,1){20}}
\put(62.67,6){\circle*{1.33}}
\put(83,6){\circle*{1.33}}
\put(83,26){\circle*{1.33}}
\put(62.67,26){\circle*{1.33}}
\put(73.33,0){\makebox(0,0)[cc]{($a$)}}
\put(113.33,0){\makebox(0,0)[cc]{(b)}}
\put(93.33,-4){\makebox(0,0)[cc]{Fig 28}}

\put(70,16){\makebox(0,0)[cc]{$v_{i}$}}
\put(103,29){\makebox(0,0)[cc]{$v_{i}$}}
\put(102.67,6){\line(1,0){20.67}}
\put(102.67,6){\line(0,1){20.67}}
\put(123,6){\line(0,1){20.67}}
\put(102.67,26){\line(1,0){20.67}}
\put(103,26){\line(1,-1){20}}
\put(102.67,6){\circle*{1.33}}
\put(123,6){\circle*{1.33}}
\put(123,26){\circle*{1.33}}
\put(102.67,26){\circle*{1.33}}
\end{picture}\\\\\\
Now we add the values of F$_{n}$ arising from the above cases and determine $x$. Substituting the value of $x$ in \\$\frac{1}{2}(a_{ii}^{(6)}- x)$ and simplifying, we get the number of $6-$cycles each of which contains a specific vertex $v_{i}$ of G.\hfill$\Box$\\
\begin{Example}
In the graph of Fig $29$ we have, F$_{1}= 5,$ F$_{2}= 60,$ F$_{3}= 60,$ F$_{4}= 60,$ F$_{5}= 120,$ F$_{6}= 60,$ F$_{7}= 60,$ F$_{8}= 80,$ F$_{9}= 360,$ F$_{10}= 360,$ F$_{11}= 180,$ F$_{12}= 120,$ F$_{13}= 240,$ F$_{14}= 120,$ F$_{15}= 240,$ F$_{16}= 240,$ F$_{17}= 120$. So, we have $x=$ $2485$. Consequently, by Theorem $3.1$, the number of $ ~ 6-$cycles each of which contains the vertex $v_{1}$ in the graph of Fig $29$ is $60$.
\end{Example}
\unitlength=0.9mm
\special{em:linewidth 0.4pt}
\linethickness{0.4pt}
\begin{picture}(95,44)
\put(82.67,2){\line(0,1){40}}
\put(103.33,2){\line(0,1){40}}
\multiput(82.67,2)(.03,.0197){1010}{\line(1,0){.0351123596}}
\multiput(113.70,22)(-.03,.0197){1010}{\line(1,0){.0351123596}}
\multiput(103.33,2)(-.03,.0587){694}{\line(1,0){.0351123596}}
\multiput(82.67,2)(.03,.0585){694}{\line(1,0){.0351123596}}
\multiput(103.33,2)(-.03,.0198){1010}{\line(1,0){.0351123596}}
\multiput(72.40,22)(.03,.0198){1010}{\line(1,0){.0351123596}}
\put(82.67,2){\line(-1,2){10}}
\put(103.33,2){\line(1,2){10}}
\put(92.90,22){\line(1,0){20.67}}
\put(113.70,22){\line(-1,2){10}}
\put(82.67,42){\line(1,0){20.67}}
\put(82.67,42){\line(-1,-2){10}}
\put(92.90,22){\line(-1,0){20.67}}
\put(82.67,2){\line(1,0){20.67}}
\put(82.67,42){\circle*{1.33}}
\put(103.33,42){\circle*{1.33}}
\put(72.80,22){\circle*{1.33}}
\put(113.4,22){\circle*{1.33}}
\put(82.67,2){\circle*{1.33}}
\put(103.33,2){\circle*{1.33}}
\put(93.33,-4){\makebox(0,0)[cc]{Fig 29}}
\put(79,42){\makebox(0,0)[cc]{$v_{1}$}}
\put(107,42){\makebox(0,0)[cc]{$v_{2}$}}
\put(69,22){\makebox(0,0)[cc]{$v_{6}$}}
\put(117,22){\makebox(0,0)[cc]{$v_{3}$}}
\put(79,2){\makebox(0,0)[cc]{$v_{5}$}}
\put(107,2){\makebox(0,0)[cc]{$v_{4}$}}
\end{picture}\\
\begin{Theorem}
If G is a simple graph with n vertices and the adjacency matrix  A= $[a_{ij}]$, then the number of $7-$cycles each of which contains a specific vertex $v_{i}$ of G is $\frac{1}{2}$ $ (a_{ii}^{(7)}-$ x$)$, where x is equal to $\displaystyle{\sum_{i=1}^{30}}$$F_{i}$ in the cases those are considered below.
\end{Theorem}
\noindent {\bf Proof:}
{The number of $7-$cycles each of which contains a specific vertex $v_{i}$ of the graph G is equal to $\frac{1}{2}$ $(a_{ii}^{(7)}-x)$, where $x$ is the number of closed walks of length $7$ form the vertex $v_{i}$ to $v_{i}$ that are not $7-$cycles. To find $x$, we have 30 cases as considered below; the cases are based on the configurations-(subgraphs) that generate $v_{i}-v_{i}$ walks of length 7 that are not cycles. In each case, N denote the number of walks of length $7$ from $v_{i}$ to $v_{i}$ that are not cycles in the corresponding subgraph, M denote the number of subgraphs of G of the same configuration and  F$_{n}$, (n= 1, 2, ...) denote the total number of $v_{i}-v_{i}$ walks of length $7$ that are not cycles in all possible subgraphs of G of the same configuration. However, in the cases with more than one figure (Cases 9, 10, ..., 18, 20,..., 30), N, M and F$_{n}$ are based on the first graph of the respective figures and P$_{1}$, P$_{2}$,... denote the number of subgraphs of G which don't have the same configuration as the first graph but are counted in M. It is clear that  F$_{n}$ is equal to N $\times$ (M$-$ P$_{1}-$P$_{2}-$...). To find N in each case, we have to include in any walk, all the edges and the vertices of the corresponding subgraphs at least once.\\
\noindent {\bf Case 1:}
{For the configuration of Fig 30, N= 42, M= $\frac{1}{2}\displaystyle{\sum_{j\neq i}}a_{ij}^{(2)}a_{ij}$ and F$_{1}=$ 21$\displaystyle{\sum_{j\neq i}}a_{ij}^{(2)}a_{ij}$.\\
\unitlength=0.9mm
\special{em:linewidth 0.4pt}
\linethickness{0.4pt}
\begin{picture}(83.33,28)
\put(82.67,2){\line(1,2){10}}
\put(103.33,2){\line(-1,2){10}}
\put(82.67,2){\line(1,0){20.67}}
\put(93,22){\circle*{1.33}}
\put(82.67,2){\circle*{1.33}}
\put(103.33,2){\circle*{1.33}}
\put(93.8,25){\makebox(0,0)[cc]{$v_{i}$}}
\put(93.33,-4){\makebox(0,0)[cc]{Fig 30}}
\end{picture}\\\\\\
\noindent {\bf Case 2:}
{For the configuration of Fig $31$, N=$ ~34$, M= $ \frac{1}{2}a_{ii}^{(3)}(d_{i}-2)$ and F$_{2}$= $ 17a_{ii}^{(3)} (d_{i}-2)$.\\\\
\unitlength=0.9mm
\special{em:linewidth 0.4pt}
\linethickness{0.4pt}
\begin{picture}(83.33,25)
\put(82.67,2){\line(1,2){10}}
\put(103.33,2){\line(-1,2){10}}
\put(92.90,22){\line(1,0){20.67}}
\put(82.67,2){\line(1,0){20.67}}
\put(92.90,22){\circle*{1.33}}
\put(113.90,22){\circle*{1.33}}
\put(82.67,2){\circle*{1.33}}
\put(103.33,2){\circle*{1.33}}
\put(93.33,-2){\makebox(0,0)[cc]{Fig 31}}
\put(92.90,24){\makebox(0,0)[cc]{$v_{i}$}}
\end{picture}\\\\\\
\noindent {\bf Case 3:}
{For the configuration of Fig 32, N= 18, M=$\displaystyle{\sum_{j\neq i}}a_{ij}^{(2)}a_{ij}(d_{j}-2)$ and F$_{3}$= 18$\displaystyle{\sum_{j\neq i}}a_{ij}^{(2)}a_{ij}(d_{j}-2)$.\\
\unitlength=0.9mm
\special{em:linewidth 0.4pt}
\linethickness{0.4pt}
\begin{picture}(83.33,25)
\put(82.67,2){\line(1,2){10}}
\put(103.33,2){\line(-1,2){10}}
\put(92.90,22){\line(1,0){20.67}}
\put(82.67,2){\line(1,0){20.67}}
\put(92.90,22){\circle*{1.33}}
\put(113.90,22){\circle*{1.33}}
\put(82.67,2){\circle*{1.33}}
\put(103.33,2){\circle*{1.33}}
\put(93.33,-2){\makebox(0,0)[cc]{Fig 32}}
\put(80,2){\makebox(0,0)[cc]{$v_{i}$}}
\end{picture}\\\\\\
\noindent {\bf Case 4:}
{For the configuration of Fig 33, N= 4, M=$\displaystyle{\sum_{j\neq i}}(\frac{1}{2}a_{jj}^{(3)}a_{ij}-a_{ij}^{(2)}a_{ij})(d_{j}-3)$   and F$_{4}$= 4$\displaystyle{\sum_{j\neq i}}(\frac{1}{2}a_{jj}^{(3)}a_{ij}-a_{ij}^{(2)}a_{ij})(d_{j}-3)$.\\
\unitlength=0.9mm
\special{em:linewidth 0.4pt}
\linethickness{0.4pt}
\begin{picture}(83.33,20)
\put(82.67,2){\line(1,2){10}}
\put(103.33,2){\line(-1,2){10}}
\put(92.90,22){\line(1,0){20.67}}
\put(92.90,22){\line(-1,0){20.67}}
\put(82.67,2){\line(1,0){20.67}}
\put(92.90,22){\circle*{1.33}}
\put(72.90,22){\circle*{1.33}}
\put(113.90,22){\circle*{1.33}}
\put(82.67,2){\circle*{1.33}}
\put(103.33,2){\circle*{1.33}}
\put(93.33,-2){\makebox(0,0)[cc]{Fig 33}}
\put(72.90,25){\makebox(0,0)[cc]{$v_{i}$}}
\end{picture}\\\\\\
\noindent {\bf Case 5:}
{For the configuration of Fig 34, N= 12, M= $ \frac{1}{2}$$a_{ii}^{(3)}\left(\begin{array}{c}d_{i}-2\\2\end{array}\right)$ and F$_{5}$= 6$a_{ii}^{(3)}\left(\begin{array}{c}d_{i}-2\\2\end{array}\right)$.\\
\unitlength=0.9mm
\special{em:linewidth 0.4pt}
\linethickness{0.4pt}
\begin{picture}(83.33,30)
\put(82.67,2){\line(1,2){10}}
\put(103.33,2){\line(-1,2){10}}
\put(92.90,22){\line(1,0){20.67}}
\put(92.90,22){\line(-1,0){20.67}}
\put(82.67,2){\line(1,0){20.67}}
\put(92.90,22){\circle*{1.33}}
\put(72.90,22){\circle*{1.33}}
\put(113.90,22){\circle*{1.33}}
\put(82.67,2){\circle*{1.33}}
\put(103.33,2){\circle*{1.33}}
\put(93.33,-2){\makebox(0,0)[cc]{Fig 34}}
\put(92.90,25){\makebox(0,0)[cc]{$v_{i}$}}
\end{picture}\\\\\\
\noindent {\bf Case 6:}
For the configuration of Fig 35, N= 4, M=$\displaystyle{\sum_{j\neq i}} a_{ij}^{(2)}a_{ij} \left(\begin{array}{c}d_{j}-2\\2\end{array}\right)$ and F$_{6}$= 4$\displaystyle{\sum_{j\neq i}} a_{ij}^{(2)}a_{ij} \left(\begin{array}{c}d_{j}-2\\2\end{array}\right)$.\\
\unitlength=0.9mm
\unitlength=0.9mm
\special{em:linewidth 0.4pt}
\linethickness{0.4pt}
\begin{picture}(83.33,25)
\put(82.67,2){\line(1,2){10}}
\put(103.33,2){\line(-1,2){10}}
\put(92.90,22){\line(1,0){20.67}}
\put(92.90,22){\line(-1,0){20.67}}
\put(82.67,2){\line(1,0){20.67}}
\put(92.90,22){\circle*{1.33}}
\put(72.90,22){\circle*{1.33}}
\put(113.90,22){\circle*{1.33}}
\put(82.67,2){\circle*{1.33}}
\put(103.33,2){\circle*{1.33}}
\put(93.33,-2){\makebox(0,0)[cc]{Fig 35}}
\put(79.90,2){\makebox(0,0)[cc]{$v_{i}$}}
\end{picture}\\\\\\
\noindent {\bf Case 7:}
For the configuration of Fig 36, N= 32, M=$\displaystyle{\sum_{j\neq i}}$$\left(\begin{array}{c}a_{ij}^{(2)}\\2\end{array}\right)$$ a_{ij}$ and F$_{7}$= 32$\displaystyle{\sum_{j\neq i}}$$\left(\begin{array}{c}a_{ij}^{(2)}\\2\end{array}\right)$$ a_{ij}$.\\
\unitlength=0.9mm
\special{em:linewidth 0.4pt}
\linethickness{0.4pt}
\begin{picture}(83,28)
\put(82.67,2){\line(1,2){10}}
\put(103.33,2){\line(-1,2){10}}
\put(103.33,2){\line(1,2){10}}
\put(92.90,22){\line(1,0){20.67}}
\put(82.67,2){\line(1,0){20.67}}
\put(92.90,22){\circle*{1.33}}
\put(113.90,22){\circle*{1.33}}
\put(82.67,2){\circle*{1.33}}
\put(103.33,2){\circle*{1.33}}
\put(93.33,-2){\makebox(0,0)[cc]{Fig 36}}
\put(92.90,25){\makebox(0,0)[cc]{$v_{i}$}}
\end{picture}\\\\
\noindent {\bf Case 8:}
For the configuration of Figure $37$, N=$24$, M= $\displaystyle{\sum_{j\neq i}} $$\left(\begin{array}{c}a_{ij}^{(2)}\\3\end{array}\right)a_{ij} $, F$_{8}$= $24$$\displaystyle{\sum_{j\neq i}}$$\left(\begin{array}{c}a_{ij}^{(2)}\\3\end{array}\right)$$ a_{ij}$.\\
\unitlength=0.9mm
\special{em:linewidth 0.4pt}
\linethickness{0.4pt}
\begin{picture}(83.33,38)
\put(82.67,2){\line(1,0){20}}
\put(102.67,2){\line(-1,1){20}}
\put(82.67,2){\line(0,1){20}}
\put(102.67,2){\line(0,1){20}}
\put(102.67,2){\line(1,3){10}}
\put(82.67,22){\line(1,0){20}}
\put(82.67,22){\line(3,1){30}}
\put(82.67,2){\circle*{1.33}}
\put(102.67,2){\circle*{1.33}}
\put(102.67,22){\circle*{1.33}}
\put(82.67,22){\circle*{1.33}}
\put(112.3,31.6){\circle*{1.33}}
\put(80,22){\makebox(0,0)[cc]{$v_{i}$}}
\put(91.90,-4){\makebox(0,0)[cc]{Fig 37}}
\end{picture}\\\\\\
\noindent {\bf Case 9:}
For the configuration of Figure $38(a)$, N= $12$, M= $\frac{1}{2}\displaystyle{\sum_{k\neq j,\ j,k\neq i}}\left(\begin{array}{c}a_{jk}^{(2)}\\2\end{array}\right) a_{jk}a_{ij}a_{ik}$. Let $P_{1}$  denote the number of all subgraphs of G that have the same configuration as the graph of Fig 38(b) and are counted in M. Thus P$_{1}$= $1\times(\frac{1}{24}$ F$_{10})$, where $\frac{1}{24}$ F$_{10}$ is the number of subgraphs of G that have the same configuration as the graph of Fig 38(b) and this subgraph is counted only once in M. Consequently, F$_{9}$$= 6$$\displaystyle{\sum_{k\neq j,\ j,k\neq i}}\left(\begin{array}{c}a_{jk}^{(2)}\\2\end{array}\right) a_{jk}a_{ij}a_{ik}- $ $\frac{1}{2}$ F$_{10}$.\\
\unitlength=0.9mm
\special{em:linewidth 0.4pt}
\linethickness{0.4pt}
\begin{picture}(83.33,32)
\put(62.67,2){\line(1,0){20}}
\put(82.67,2){\line(-1,1){20}}
\put(62.67,2){\line(0,1){20}}
\put(82.67,2){\line(0,1){20}}
\put(82.67,2){\line(1,3){10}}
\put(62.67,22){\line(1,0){20}}
\put(62.67,22){\line(3,1){30}}
\put(62.67,2){\circle*{1.33}}
\put(82.67,2){\circle*{1.33}}
\put(82.67,22){\circle*{1.33}}
\put(62.67,22){\circle*{1.33}}
\put(92.3,31.6){\circle*{1.33}}
\put(92,34){\makebox(0,0)[cc]{$v_{i}$}}
\put(91.90,-4){\makebox(0,0)[cc]{Fig 38}}

\put(102.67,2){\line(1,2){10}}
\put(123.33,2){\line(-1,2){10}}
\put(123.33,2){\line(1,2){10}}
\put(112.90,22){\line(1,0){20.67}}
\put(102.67,2){\line(1,0){20.67}}
\put(112.90,22){\circle*{1.33}}
\put(133.90,22){\circle*{1.33}}
\put(102.67,2){\circle*{1.33}}
\put(123.33,2){\circle*{1.33}}
\put(132.90,25){\makebox(0,0)[cc]{$v_{i}$}}
\put(73.33,-2){\makebox(0,0)[cc]{($a$)}}
\put(113.33,-2){\makebox(0,0)[cc]{(b)}}
\end{picture}\\\\\\
\noindent {\bf Case 10:}
For the configuration of Fig 39(a), N= 24, M= $\frac{1}{2}\displaystyle{\sum_{k\neq j,\ j,k\neq i}}a_{jk}^{(2)}a_{ij}a_{ik}a_{jk}$. Let $P_{1}$  denote the number of all subgraphs of G that have the same configuration as the graph of Fig 39(b) and are counted in M. Thus P$_{1}$= $1\times(\frac{1}{42}$ F$_{1})$, where $\frac{1}{42}$ F$_{1}$ is the number of subgraphs of G that have the same configuration as the graph of Fig 39(b) and this subgraph is counted only once in M. Consequently, F$_{10}$$= 12 \displaystyle{\sum_{k\neq j,\ j,k\neq i}}a_{jk}^{(2)}a_{ij}a_{ik}a_{jk}- $ $\frac{4}{7}$ F$_{1}$ .\\
\unitlength=0.9mm
\special{em:linewidth 0.4pt}
\linethickness{0.4pt}
\begin{picture}(83,30)
\put(62.67,2){\line(1,2){10}}
\put(83.33,2){\line(-1,2){10}}
\put(83.33,2){\line(1,2){10}}
\put(72.90,22){\line(1,0){20.67}}
\put(62.67,2){\line(1,0){20.67}}
\put(72.90,22){\circle*{1.33}}
\put(93.90,22){\circle*{1.33}}
\put(62.67,2){\circle*{1.33}}
\put(83.33,2){\circle*{1.33}}
\put(93.33,-4){\makebox(0,0)[cc]{Fig 39}}
\put(92.90,25){\makebox(0,0)[cc]{$v_{i}$}}

\put(102.67,2){\line(1,2){10}}
\put(123.33,2){\line(-1,2){10}}
\put(102.67,2){\line(1,0){20.67}}
\put(113,22){\circle*{1.33}}
\put(102.67,2){\circle*{1.33}}
\put(123.33,2){\circle*{1.33}}
\put(113.5,25){\makebox(0,0)[cc]{$v_{i}$}}
\put(73.33,-2){\makebox(0,0)[cc]{($a$)}}
\put(113.33,-2){\makebox(0,0)[cc]{(b)}}
\end{picture}\\\\\\
\noindent {\bf Case 11:}
{For the configuration of Fig 40(a), N= 14, M= $\frac{1}{2}\displaystyle{\sum_{j\neq i}}a_{jj}^{(3)}a_{ij}$. Let P$_{1}$ denote the number of all subgraphs of G that have the same configuration as the graph of Fig 40(b) and are counted in M. Thus P$_{1}$= 2$\times (\frac{1}{42}$F$_{1})$, where $\frac{1}{42}$ F$_{1}$ is the number of subgraphs of G that have the same configuration as the graph of Fig 40(b) and 2 is the number of times that this subgraph is counted in M. Consequently, F$_{11}$= 7$\displaystyle{\sum_{j\neq i}}a_{jj}^{(3)}a_{ij}-$ $\frac{2}{3}$ F$_{1}$.\\
\unitlength=0.9mm
\special{em:linewidth 0.4pt}
\linethickness{0.4pt}
\begin{picture}(83.33,28)
\put(62.67,2){\line(1,2){10}}
\put(83.33,2){\line(-1,2){10}}
\put(72.90,22){\line(1,0){20.67}}
\put(62.67,2){\line(1,0){20.67}}
\put(72.90,22){\circle*{1.33}}
\put(93.90,22){\circle*{1.33}}
\put(62.67,2){\circle*{1.33}}
\put(83.33,2){\circle*{1.33}}
\put(73.33,-2){\makebox(0,0)[cc]{($a$)}}
\put(93.90,25){\makebox(0,0)[cc]{$v_{i}$}}
\put(102.67,2){\line(1,2){10}}
\put(123.33,2){\line(-1,2){10}}
\put(102.67,2){\line(1,0){20.67}}
\put(112.90,22){\circle*{1.33}}
\put(102.67,2){\circle*{1.33}}
\put(123.33,2){\circle*{1.33}}
\put(113.33,-2){\makebox(0,0)[cc]{(b)}}
\put(112.90,25){\makebox(0,0)[cc]{$v_{i}$}}
\put(94,-4){\makebox(0,0)[cc]{Fig 40}}
\end{picture}\\\\\\
\noindent {\bf Case 12:}
{For the configuration of Fig $41(a)$, N= 2, M= $\frac{1}{2}\displaystyle{\sum_{j\neq i}}a_{ij}^{(2)}a_{jj}^{(3)}$. Let P$_{1}$ denote the number of all subgraphs of G that have the same configuration as the graph of Fig 41(b) and are counted in M. Thus P$_{1}=2\times(\frac{1}{42}$F$_{1})$, where $\frac{1}{42}$ F$_{1}$ is the number of subgraphs of G that have the same configuration as the graph of Fig 41(b) and 2 is the number of times that this subgraph is counted in M. Let P$_{2}$ denote the number of all subgraphs of G that have the same configuration as the graph of Fig 41(c) and are counted in M. Thus P$_{2}=2\times( \frac{1}{14}$ F$_{11})$, where $\frac{1}{14}$ F$_{11}$ is the number of subgraphs of G that have the same configuration as the graph of Fig 41(c) and 2 is the number of times that this subgraph is counted in M. Let P$_{3}$ denote the number of all subgraphs of G that have the same configuration as the graph of Fig 41(d) and are counted in M. Thus P$_{3}=2\times(\frac{1}{32}$ F$_{7})$, where $\frac{1}{32}$ F$_{7}$ is the number of subgraphs of G that have the same configuration as the graph of Fig 41(d) and 2 is the number of times that this subgraph is counted in M. Consequently, F$_{12}$= $\displaystyle{\sum_{j\neq i}}a_{ij}^{(2)}a_{jj}^{(3)}$$ -$  $\frac{2}{21}$ F$_{1}- \frac{2}{7}$ F$_{11}-\frac{1}{8}$ F$_{7}$.\\
\unitlength=0.9mm
\special{em:linewidth 0.4pt}
\linethickness{0.4pt}
\begin{picture}(136,30)
\put(118.33,2){\line(-1,2){10}}
\put(107.90,22){\line(1,0){20.67}}
\put(27,22){\line(-1,-2){10}}
\put(97.90,2){\line(1,2){10}}
\put(47,22){\line(-1,0){20.67}}
\put(97.67,2){\line(1,0){20.67}}
\put(26.90,22){\circle*{1.33}}
\put(17,2){\circle*{1.33}}
\put(97.90,2){\circle*{1.33}}
\put(107.90,22){\circle*{1.33}}
\put(117.67,2){\circle*{1.33}}
\put(128.33,22){\circle*{1.33}}
\put(92.33,-6){\makebox(0,0)[cc]{Fig 41}}
\put(107.67,-2){\makebox(0,0)[cc]{$( c)$}}
\put(14,2){\makebox(0,0)[cc]{$v_{i}$}}
\put(128,25){\makebox(0,0)[cc]{$v_{i}$}}
\put(77.90,25){\makebox(0,0)[cc]{$v_{i}$}}
\put(143.90,24){\makebox(0,0)[cc]{$v_{i}$}}
\put(37.67,2){\line(1,2){10}}
\put(58.33,2){\line(-1,2){10}}
\put(37.67,2){\line(1,0){20.67}}
\put(47.90,22){\circle*{1.33}}
\put(37.67,2){\circle*{1.33}}
\put(58.33,2){\circle*{1.33}}
\put(47,-2){\makebox(0,0)[cc]{$( a )$}}
\unitlength=0.9mm
\special{em:linewidth 0.4pt}
\linethickness{0.4pt}
\put(132.67,2){\line(1,2){10}}
\put(153.33,2){\line(-1,2){10}}
\put(153.33,2){\line(1,2){10}}
\put(142.90,22){\line(1,0){20.67}}
\put(132.67,2){\line(1,0){20.67}}
\put(142.90,22){\circle*{1.33}}
\put(163.90,22){\circle*{1.33}}
\put(132.67,2){\circle*{1.33}}
\put(153.33,2){\circle*{1.33}}
\put(143,-2){\makebox(0,0)[cc]{$( d)$}}
\put(67.67,2){\line(1,2){10}}
\put(88.33,2){\line(-1,2){10}}
\put(67.67,2){\line(1,0){20.67}}
\put(77.90,22){\circle*{1.33}}
\put(67.67,2){\circle*{1.33}}
\put(88.33,2){\circle*{1.33}}
\put(77.33,-2){\makebox(0,0)[cc]{$(b)$}}
\end{picture}\\\\\\\\
\noindent {\bf Case 13:}
{For the configuration of Fig $42(a)$, N= 4, M= $\frac{1}{2}\displaystyle{\sum_{j\neq i}}a_{jj}^{(3)}a_{ij}(d_{i}-1)$. Let P$_{1}$ denote the number of all subgraphs of G that have the same configuration as the graph of Fig 42(b) and are counted in M. Thus P$_{1}=$ $2\times(\frac{1}{42}$F$_{1})$, where $\frac{1}{42}$ F$_{1}$ is the number of subgraphs of G that have the same configuration as the graph of Fig 42(b) and 2 is the number of times that this subgraph is counted in M. Let P$_{2}$ denote the number of all subgraphs of G that have the same configuration as the graph of Fig 42(c) and are counted in M. Thus P$_{2}=$ $2\times(\frac{1}{34}$ F$_{2})$, where $\frac{1}{34}$ F$_{2}$ is the number of subgraphs of G that have the same configuration as the graph of Fig 42(c) and 2 is the number of times that this subgraph is counted in M. Let P$_{3}$ denote the number of all subgraphs of G that have the same configuration as the graph of Fig 42(d) and are counted in M. Thus P$_{3}=$ $2\times(\frac{1}{24}$ F$_{10})$, where $\frac{1}{24}$ F$_{10}$ is the number of subgraphs of G that have the same configuration as the graph of Fig 42(d) and 2 is the number of times that this subgraph is counted in M. Consequently, F$_{13}$= $2\displaystyle{\sum_{j\neq i}}a_{jj}^{(3)}a_{ij}(d_{i}-1)-$ $\frac{4}{21}$ F$_{1}- \frac{4}{17}$ F$_{2}- \frac{1}{3}$ F$_{10}$.\\
\unitlength=0.9mm
\special{em:linewidth 0.4pt}
\linethickness{0.4pt}
\begin{picture}(136,30)
\put(118.33,2){\line(-1,2){10}}
\put(107.90,22){\line(1,0){20.67}}
\put(27,22){\line(-1,-2){10}}
\put(97.90,2){\line(1,2){10}}
\put(47,22){\line(-1,0){20.67}}
\put(97.67,2){\line(1,0){20.67}}
\put(26.90,22){\circle*{1.33}}
\put(17,2){\circle*{1.33}}
\put(97.90,2){\circle*{1.33}}
\put(107.90,22){\circle*{1.33}}
\put(117.67,2){\circle*{1.33}}
\put(128.33,22){\circle*{1.33}}
\put(92.33,-6){\makebox(0,0)[cc]{Fig 42}}
\put(107.67,-2){\makebox(0,0)[cc]{$( c)$}}
\put(78.50,24){\makebox(0,0)[cc]{$v_{i}$}}
\put(108,24){\makebox(0,0)[cc]{$v_{i}$}}
\put(27.90,24){\makebox(0,0)[cc]{$v_{i}$}}
\put(163.90,24){\makebox(0,0)[cc]{$v_{i}$}}
\put(37.67,2){\line(1,2){10}}
\put(58.33,2){\line(-1,2){10}}
\put(37.67,2){\line(1,0){20.67}}
\put(47.90,22){\circle*{1.33}}
\put(37.67,2){\circle*{1.33}}
\put(58.33,2){\circle*{1.33}}
\put(47,-2){\makebox(0,0)[cc]{$( a )$}}
\unitlength=0.9mm
\special{em:linewidth 0.4pt}
\linethickness{0.4pt}
\put(132.67,2){\line(1,2){10}}
\put(153.33,2){\line(-1,2){10}}
\put(153.33,2){\line(1,2){10}}
\put(142.90,22){\line(1,0){20.67}}
\put(132.67,2){\line(1,0){20.67}}
\put(142.90,22){\circle*{1.33}}
\put(163.90,22){\circle*{1.33}}
\put(132.67,2){\circle*{1.33}}
\put(153.33,2){\circle*{1.33}}
\put(143,-2){\makebox(0,0)[cc]{$( d)$}}
\put(67.67,2){\line(1,2){10}}
\put(88.33,2){\line(-1,2){10}}
\put(67.67,2){\line(1,0){20.67}}
\put(77.90,22){\circle*{1.33}}
\put(67.67,2){\circle*{1.33}}
\put(88.33,2){\circle*{1.33}}
\put(77.33,-2){\makebox(0,0)[cc]{$(b)$}}
\end{picture}\\\\\\\\
\noindent {\bf Case 14:}
{For the configuration of Fig $43(a)$, N= 4, M= $\frac{1}{2}\displaystyle{\sum_{j\neq i}}a_{ii}^{(3)}a_{ij}^{(2)}$. Let P$_{1}$ denote the number of all subgraphs of G that have the same configuration as the graph of Fig 43(b) and are counted in M. Thus P$_{1}=$ $2\times(\frac{1}{42}$F$_{1})$, where $\frac{1}{42}$ F$_{1}$ is the number of subgraphs of G that have the same configuration as the graph of Fig 43(b) and 2 is the number of times that this subgraph is counted in M. Let P$_{2}$ denote the number of all subgraphs of G that have the same configuration as the graph of Fig 43(c) and are counted in M. Thus P$_{2}=$ $1\times(\frac{1}{18}$ F$_{3})$, where $\frac{1}{18}$ F$_{3}$ is the number of subgraphs of G that have the same configuration as the graph of Fig 43(c) and this subgraph is counted only once in M. Let P$_{3}$ denote the number of all subgraphs of G that have the same configuration as the graph of Fig 43(d) and are counted in M. Thus P$_{3}=$ $2\times(\frac{1}{32}$ F$_{7})$, where $\frac{1}{32}$ F$_{7}$ is the number of subgraphs of G that have the same configuration as the graph of Fig 43(d) and 2 is the number of times that this subgraph is counted in M. \\Consequently, F$_{14}$= $2\displaystyle{\sum_{j\neq i}}a_{ii}^{(3)}a_{ij}^{(2)}-$ $\frac{4}{21}$ F$_{1}- $ $\frac{2}{9}$ F$_{3}-$$ \frac{1}{4}$ F$_{7}$.\\
\unitlength=0.9mm
\special{em:linewidth 0.4pt}
\linethickness{0.4pt}
\begin{picture}(136,30)
\put(118.33,2){\line(-1,2){10}}
\put(107.90,22){\line(1,0){20.67}}
\put(27,22){\line(-1,-2){10}}
\put(97.90,2){\line(1,2){10}}
\put(47,22){\line(-1,0){20.67}}
\put(97.67,2){\line(1,0){20.67}}
\put(26.90,22){\circle*{1.33}}
\put(17,2){\circle*{1.33}}
\put(97.90,2){\circle*{1.33}}
\put(107.90,22){\circle*{1.33}}
\put(117.67,2){\circle*{1.33}}
\put(128.33,22){\circle*{1.33}}
\put(92.33,-6){\makebox(0,0)[cc]{Fig 43}}
\put(107.67,-2){\makebox(0,0)[cc]{$( c)$}}
\put(78.50,24){\makebox(0,0)[cc]{$v_{i}$}}
\put(121,2){\makebox(0,0)[cc]{$v_{i}$}}
\put(47.90,24){\makebox(0,0)[cc]{$v_{i}$}}
\put(143.90,24){\makebox(0,0)[cc]{$v_{i}$}}
\put(37.67,2){\line(1,2){10}}
\put(58.33,2){\line(-1,2){10}}
\put(37.67,2){\line(1,0){20.67}}
\put(47.90,22){\circle*{1.33}}
\put(37.67,2){\circle*{1.33}}
\put(58.33,2){\circle*{1.33}}
\put(47,-2){\makebox(0,0)[cc]{$( a )$}}
\unitlength=0.9mm
\special{em:linewidth 0.4pt}
\linethickness{0.4pt}
\put(132.67,2){\line(1,2){10}}
\put(153.33,2){\line(-1,2){10}}
\put(153.33,2){\line(1,2){10}}
\put(142.90,22){\line(1,0){20.67}}
\put(132.67,2){\line(1,0){20.67}}
\put(142.90,22){\circle*{1.33}}
\put(163.90,22){\circle*{1.33}}
\put(132.67,2){\circle*{1.33}}
\put(153.33,2){\circle*{1.33}}
\put(143,-2){\makebox(0,0)[cc]{$( d)$}}
\put(67.67,2){\line(1,2){10}}
\put(88.33,2){\line(-1,2){10}}
\put(67.67,2){\line(1,0){20.67}}
\put(77.90,22){\circle*{1.33}}
\put(67.67,2){\circle*{1.33}}
\put(88.33,2){\circle*{1.33}}
\put(77.33,-2){\makebox(0,0)[cc]{$(b)$}}
\end{picture}\\\\\\\\
\noindent {\bf Case 15:}
{For the configuration of Fig $44(a)$, N= 2, M=$\displaystyle{\sum_{j\neq k,\ j,k\neq i}}a_{ij}^{(2)}a_{jk}^{(2)}a_{ij}$. Let P$_{1}$ denote the number of all subgraphs of G that have the same configuration as the graph of Fig 44(b) and are counted in M. Thus P$_{1}=2\times(\frac{1}{42}$F$_{1})$, where $\frac{1}{42}$ F$_{1}$ is the number of subgraphs of G that have the same configuration as the graph of Fig 44(b) and 2 is the number of times that this subgraph is counted in M. Let P$_{2}$ denote the number of all subgraphs of G that have the same configuration as the graph of Fig 44(c) and are counted in M. Thus P$_{2}=$ $2\times(\frac{1}{34}$ F$_{2})$, where $\frac{1}{34}$ F$_{2}$ is the number of subgraphs of G that have the same configuration as the graph of Fig 44(c) and 2 is the number of times that this subgraph is counted in M. Let P$_{3}$ denote the number of all subgraphs of G that have the same configuration as the graph of Fig 44(d) and are counted in M. Thus P$_{3}=$ $2\times(\frac{1}{24}$ F$_{10})$, where $\frac{1}{24}$ F$_{10}$ is the number of subgraphs of G that have the same configuration as the graph of Fig 44(d) and 2 is the number of times that this subgraph is counted in M. 
Let P$_{4}$ denote the number of all subgraphs of G that have the same configuration as the graph of Fig 44(e) and are counted in M. Thus P$_{4}=$ $1\times(\frac{1}{18}$ F$_{3})$, where $\frac{1}{18}$ F$_{3}$ is the number of subgraphs of G that have the same configuration as the graph of Fig 44(e) and 1 is the number of times that this subgraph is counted in M.  Consequently, F$_{15}$= $2\displaystyle{\sum_{j\neq k,\ j,k\neq i}}a_{ij}^{(2)}a_{jk}^{(2)}a_{ij}$$- ~\frac{2}{21}$ F$_{1}- \frac{2}{17}$ F$_{2}- \frac{1}{6}$ F$_{10}- \frac{1}{9}$ F$_{3}$.\\
\unitlength=0.9mm
\special{em:linewidth 0.4pt}
\linethickness{0.4pt}
\begin{picture}(136,30)
\put(98.33,2){\line(-1,2){10}}
\put(87.90,22){\line(1,0){20.67}}
\put(7,22){\line(-1,-2){10}}
\put(77.90,2){\line(1,2){10}}
\put(27,22){\line(-1,0){20.67}}
\put(77.67,2){\line(1,0){20.67}}
\put(6.90,22){\circle*{1.33}}
\put(-3,2){\circle*{1.33}}
\put(77.90,2){\circle*{1.33}}
\put(87.90,22){\circle*{1.33}}
\put(97.67,2){\circle*{1.33}}
\put(108.33,22){\circle*{1.33}}
\put(88.33,-8){\makebox(0,0)[cc]{Fig 44}}
\put(87.67,-2){\makebox(0,0)[cc]{$( c)$}}
\put(41,2){\makebox(0,0)[cc]{$v_{i}$}}
\put(88,25){\makebox(0,0)[cc]{$v_{i}$}}
\put(57.90,25){\makebox(0,0)[cc]{$v_{i}$}}
\put(143.90,24){\makebox(0,0)[cc]{$v_{i}$}}
\put(17.67,2){\line(1,2){10}}
\put(38.33,2){\line(-1,2){10}}
\put(17.67,2){\line(1,0){20.67}}
\put(27.90,22){\circle*{1.33}}
\put(17.67,2){\circle*{1.33}}
\put(38.33,2){\circle*{1.33}}
\put(27,-2){\makebox(0,0)[cc]{$( a )$}}
\unitlength=0.9mm
\special{em:linewidth 0.4pt}
\linethickness{0.4pt}
\put(112.67,2){\line(1,2){10}}
\put(133.33,2){\line(-1,2){10}}
\put(133.33,2){\line(1,2){10}}
\put(122.90,22){\line(1,0){20.67}}
\put(112.67,2){\line(1,0){20.67}}
\put(122.90,22){\circle*{1.33}}
\put(143.90,22){\circle*{1.33}}
\put(112.67,2){\circle*{1.33}}
\put(133.33,2){\circle*{1.33}}
\put(123,-2){\makebox(0,0)[cc]{$( d)$}}
\put(47.67,2){\line(1,2){10}}
\put(68.33,2){\line(-1,2){10}}
\put(47.67,2){\line(1,0){20.67}}
\put(57.90,22){\circle*{1.33}}
\put(47.67,2){\circle*{1.33}}
\put(68.33,2){\circle*{1.33}}
\put(57.33,-2){\makebox(0,0)[cc]{$(b)$}}

\put(152.67,2){\line(1,2){10}}
\put(173.33,2){\line(-1,2){10}}
\put(162.90,22){\line(1,0){20.67}}
\put(152.67,2){\line(1,0){20.67}}
\put(162.90,22){\circle*{1.33}}
\put(183.90,22){\circle*{1.33}}
\put(152.67,2){\circle*{1.33}}
\put(173.33,2){\circle*{1.33}}
\put(150,2){\makebox(0,0)[cc]{$v_{i}$}}
\put(162.33,-2){\makebox(0,0)[cc]{$(e)$}}
\end{picture}\\\\\\\\
\noindent {\bf Case 16:}
For the configuration of Figure 45(a), N=$~2$, M=$\displaystyle{\sum_{j\neq k,\ j,k\neq i}}a_{jk}^{(2)} a_{jk}a_{ij}(d_{k}-2)$. Let P$_{1}$ denote the number of all subgraphs of G that have the same configuration as the graph of Figure 45(b) and are counted in M.  Thus P$_{1}=$ $2\times(\frac{1}{24}$ F$_{10})$, where $\frac{1}{24}$ F$_{10}$ is the number of subgraphs of G that have the same configuration as the graph of Fig 45(b) and 2 is the number of times that this subgraph is counted in M. Let P$_{2}$ denote the number of all subgraphs of G that have the same configuration as the graph of Fig 45(c) and are counted in M. Thus P$_{2}=$ $1\times(\frac{1}{18}$ F$_{3})$, where $\frac{1}{18}$ F$_{3}$ is the number of subgraphs of G that have the same configuration as the graph of Fig 45(c) and 1 is the number of times that this subgraph is counted in M. \\Consequently, F$_{16}$= $2\displaystyle{\sum_{j\neq k,\ j,k\neq i}}a_{jk}^{(2)} a_{jk}a_{ij}(d_{k}-2)-$ $\frac{1}{6}$ F$_{10}- \frac{1}{9}$ F$_{3}$.\\
\unitlength=0.9mm
\special{em:linewidth 0.4pt}
\linethickness{0.4pt}
\begin{picture}(100,28)
\put(103.33,2){\line(1,2){10}}
\put(82.67,2){\line(1,2){10}}
\put(103.33,2){\line(-1,2){10}}
\put(92.90,22){\line(1,0){20.67}}
\put(82.67,2){\line(1,0){20.67}}
\put(92.90,22){\circle*{1.33}}
\put(82.67,2){\circle*{1.33}}
\put(103.33,2){\circle*{1.33}}
\put(113.33,22){\circle*{1.33}}
\put(92.33,-8){\makebox(0,0)[cc]{Fig 45}}
\put(92.67,-2){\makebox(0,0)[cc]{(b)}}
\put(116,22){\makebox(0,0)[cc]{$v_{i}$}}
\put(40,22){\makebox(0,0)[cc]{$v_{i}$}}
\put(52.67,2){\line(1,2){10}}
\put(73.33,2){\line(-1,2){10}}
\put(52.67,2){\line(-1,0){20.67}}
\put(62.90,22){\line(-1,0){20.67}}
\put(52.67,2){\line(1,0){20.67}}
\put(62.90,22){\circle*{1.33}}
\put(42.90,22){\circle*{1.33}}
\put(32.90,2){\circle*{1.33}}
\put(52.67,2){\circle*{1.33}}
\put(73.33,2){\circle*{1.33}}
\put(63,-2){\makebox(0,0)[cc]{($a$)}}

\put(117.67,2){\line(1,2){10}}
\put(138.33,2){\line(-1,2){10}}
\put(127.90,22){\line(1,0){20.67}}
\put(117.67,2){\line(1,0){20.67}}
\put(127.90,22){\circle*{1.33}}
\put(148.90,22){\circle*{1.33}}
\put(117.67,2){\circle*{1.33}}
\put(138.33,2){\circle*{1.33}}
\put(115,2){\makebox(0,0)[cc]{$v_{i}$}}
\put(127.33,-2){\makebox(0,0)[cc]{$(c)$}}

\end{picture}\\\\\\\\
\noindent {\bf Case 17:}
For the configuration of Figure 46(a), N=$~4$, M= $\displaystyle{\sum_{j\neq i}}a_{ij}^{(2)} a_{ij}(d_{j}-2)(d_{i}-2)$. Let P$_{1}$ denote the number of all subgraphs of G that have the same configuration as the graph of Figure 46(b) and are counted in M. Thus P$_{1}$=  $2\times(\frac{1}{32}$ F$_{7})$, where $\frac{1}{32}$ F$_{7}$ is the number of subgraphs of G that have the same configuration as the graph of Figure 46(b) and $2$  is the number of times that this subgraph is counted in M. Consequently, F$_{17}$= $4$ $\displaystyle{\sum_{j\neq i}}$ $a_{ij}^{(2)} a_{ij}(d_{j}-2)(d_{i}-2)$ $- \frac{1}{4}$ F$_{7}$.\\
\unitlength=0.9mm
\special{em:linewidth 0.4pt}
\linethickness{0.4pt}
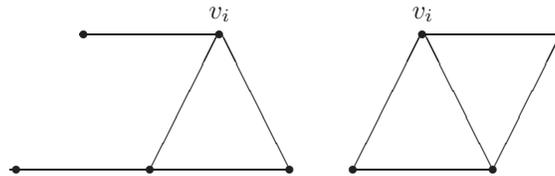
\begin{picture}(100,20)
\put(123.33,2){\line(1,2){10}}
\put(102.67,2){\line(1,2){10}}
\put(123.33,2){\line(-1,2){10}}
\put(112.90,22){\line(1,0){20.67}}
\put(102.67,2){\line(1,0){20.67}}
\put(112.90,22){\circle*{1.33}}
\put(102.67,2){\circle*{1.33}}
\put(123.33,2){\circle*{1.33}}
\put(133.33,22){\circle*{1.33}}
\put(98.33,-8){\makebox(0,0)[cc]{Fig 46}}
\put(112.67,-2){\makebox(0,0)[cc]{(b)}}

\put(72.67,2){\line(1,2){10}}
\put(93.33,2){\line(-1,2){10}}
\put(72.67,2){\line(-1,0){20.67}}
\put(82.90,22){\line(-1,0){20.67}}
\put(72.67,2){\line(1,0){20.67}}
\put(82.90,22){\circle*{1.33}}
\put(62.90,22){\circle*{1.33}}
\put(52.90,2){\circle*{1.33}}
\put(72.67,2){\circle*{1.33}}
\put(93.33,2){\circle*{1.33}}
\put(83,-2){\makebox(0,0)[cc]{($a$)}}
\put(113,25){\makebox(0,0)[cc]{$v_{i}$}}
\put(83,25){\makebox(0,0)[cc]{$v_{i}$}}

\end{picture}\\
\noindent {\bf Case 18:}
For the configuration of Figure 47(a), N=$~2$, M= $ \frac{1}{2}$ $\displaystyle{\sum_{k\neq j,\ j,k\neq i}}(d_{j}-2)(d_{k}-2)a_{ij}a_{ik}a_{jk}$. Let P$_{1}$ denotes the number of all subgraphs of G that have the same configuration as the graph of Figure 47(b) and are counted in M. Thus P$_{1}$=  $1\times(\frac{1}{24}$ F$_{10})$, where $\frac{1}{24}$ F$_{10}$ is the number of subgraphs of G that have the same configuration as the graph of Figure 47(b) and $1$  is the number of times that this subgraph is counted in M. \\Consequently, F$_{18}$= $\displaystyle{\sum_{k\neq j,\ j,k\neq i}}(d_{j}-2)(d_{k}-2)a_{ij}a_{ik}a_{jk}-  $ $\frac{1}{12}$ F$_{10}$.\\
\unitlength=0.9mm
\special{em:linewidth 0.4pt}
\linethickness{0.4pt}
\begin{picture}(100,30)
\put(123.33,2){\line(1,2){10}}
\put(102.67,2){\line(1,2){10}}
\put(123.33,2){\line(-1,2){10}}
\put(112.90,22){\line(1,0){20.67}}
\put(102.67,2){\line(1,0){20.67}}
\put(112.90,22){\circle*{1.33}}
\put(102.67,2){\circle*{1.33}}
\put(123.33,2){\circle*{1.33}}
\put(133.33,22){\circle*{1.33}}
\put(98.33,-8){\makebox(0,0)[cc]{Fig 47}}
\put(112.67,-2){\makebox(0,0)[cc]{(b)}}

\put(72.67,2){\line(1,2){10}}
\put(93.33,2){\line(-1,2){10}}
\put(72.67,2){\line(-1,0){20.67}}
\put(82.90,22){\line(-1,0){20.67}}
\put(72.67,2){\line(1,0){20.67}}
\put(82.90,22){\circle*{1.33}}
\put(62.90,22){\circle*{1.33}}
\put(52.90,2){\circle*{1.33}}
\put(72.67,2){\circle*{1.33}}
\put(93.33,2){\circle*{1.33}}
\put(83,-2){\makebox(0,0)[cc]{($a$)}}
\put(96,2){\makebox(0,0)[cc]{$v_{i}$}}
\put(133,25){\makebox(0,0)[cc]{$v_{i}$}}

\end{picture}\\\\\\\\
\noindent {\bf Case 19:}
For the configuration of Figure $48$, N= 14, M=$\frac{1}{2}$ $ [a_{ii}^{(5)} - 5 a_{ii}^{(3)} - 2 ( d_{i} - 2 ) a_{ii}^{(3)}-2$ $\displaystyle{\sum_{ j=1, j\neq i }^n}$ $a_{ij}^{(2)} a_{ij}(d_{j}-2) - 2$ $\displaystyle{\sum_{ j=1, j\neq i}^n}$ $a_{ij}(\frac{1}{2} a_{jj}^{(3)}-a_{ij}a_{ij}^{(2)})]$ (See Theorem 1.11) and F$_{19}$= 7 $ [a_{ii}^{(5)} - 5 a_{ii}^{(3)} - 2 ( d_{i} - 2 ) a_{ii}^{(3)}-2$ $\displaystyle{\sum_{ j=1, j\neq i }^n}$ $a_{ij}^{(2)} a_{ij}(d_{j}-2) - 2$ $\displaystyle{\sum_{ j=1, j\neq i}^n}$ $a_{ij}(\frac{1}{2} a_{jj}^{(3)}-a_{ij}a_{ij}^{(2)})]$.\\
\unitlength=0.9mm
\special{em:linewidth 0.4pt}
\linethickness{0.4pt}
\begin{picture}(83.33,26)
\put(85.67,2){\line(-1,2){7}}
\put(103.33,2){\line(1,2){7}}
\put(110.5,16.5){\line(-1,1){16}}
\put(78.90,16.5){\line(1,1){16}}
\put(85.67,2){\line(1,0){18}}
\put(78.90,16){\circle*{1.33}}
\put(94.7,32){\circle*{1.33}}
\put(110.5,16.5){\circle*{1.33}}
\put(85.67,2){\circle*{1.33}}
\put(103.33,2){\circle*{1.33}}
\put(93.33,-2){\makebox(0,0)[cc]{Fig 48}}
\put(95.47,35){\makebox(0,0)[cc]{$v_{i}$}}
\end{picture}\\\\
\noindent {\bf Case 20:}
For the configuration of Fig 49(a), N= $6$, M= $\frac{1}{2}$ $\displaystyle{\sum_{k\neq j,\ j,k\neq i}}$$a_{jk}^{(2)}(d_{k}-a_{jk}-1)(a_{jk} a_{ij}a_{ik})$ (See Theorem 1.5). Let P$_{1}$ denote the number of all subgraphs of G that have the same configuration as the graph of Figure 49(b) and are counted in M. Thus P$_{1}$= $2\times(\frac{1}{32}$ F$_{7})$, where $\frac{1}{32}$ F$_{7}$ is the number of subgraphs of G that have the same configuration as the graph of Figure 49(b) and $2$ is the number of times that this subgraph is counted in M. Consequently, F$_{20}$= $3\displaystyle{\sum_{k\neq j,\ j,k\neq i}}a_{jk}^{(2)}(d_{k}-a_{jk}-1)(a_{jk} a_{ij}a_{ik})-$ $\frac{3}{8}$ F$_{7}$.\\
\unitlength=0.9mm
\special{em:linewidth 0.4pt}
\linethickness{0.4pt}
\begin{picture}(83.33,35)
\put(70.67,2){\line(-1,2){7}}
\put(88.33,2){\line(1,2){7}}
\put(95.5,16.5){\line(-1,1){16}}
\put(63.90,16.5){\line(1,1){16}}
\put(70.67,2){\line(1,0){18}}
\put(63.90,16){\line(1,0){31}}
\put(63.90,16){\circle*{1.33}}
\put(79.7,32){\circle*{1.33}}
\put(95.5,16.5){\circle*{1.33}}
\put(70.67,2){\circle*{1.33}}
\put(88.33,2){\circle*{1.33}}
\put(80.47,35){\makebox(0,0)[cc]{$v_{i}$}}
\put(106.97,24){\makebox(0,0)[cc]{$v_{i}$}}
\put(92.33,-6){\makebox(0,0)[cc]{Fig 49}}
\put(107.67,-2){\makebox(0,0)[cc]{(b)}}
\put(79.67,-2){\makebox(0,0)[cc]{($a$)}}
\put(118.33,2){\line(1,2){10}}
\put(97.67,2){\line(1,2){10}}
\put(118.33,2){\line(-1,2){10}}
\put(107.90,22){\line(1,0){20.67}}
\put(97.67,2){\line(1,0){20.67}}
\put(107.90,22){\circle*{1.33}}
\put(97.67,2){\circle*{1.33}}
\put(118.33,2){\circle*{1.33}}
\put(128.33,22){\circle*{1.33}}
\end{picture}\\\\\\\\
\noindent {\bf Case 21:}
For the configuration of Fig 50(a), N= $12$, M=$\displaystyle{\sum_{j\neq i}}$$a_{ij}^{(2)}(d_{j}-a_{ij}-1)(a_{ij}^{(2)}a_{ij})$ (See Theorem 1.7). Let P$_{1}$ denote the number of all subgraphs of G that have the same configuration as the graph of Figure 50(b) and are counted in M. Thus P$_{1}$= $2\times(\frac{1}{32}$ F$_{7})$, where $\frac{1}{32}$ F$_{7}$ is the number of subgraphs of G that have the same configuration as the graph of Figure 50(b) and $2$ is the number of times that this subgraph is counted in M. Let P$_{2}$ denote the number of all subgraphs of G that have the same configuration as the graph of Figure 50(c) and are counted in M. Thus P$_{2}$= $2\times(\frac{1}{24}$ F$_{10})$, where $\frac{1}{24}$ F$_{10}$ is the number of subgraphs of G that have the same configuration as the graph of Figure 50(c) and $2$ is the number of times that this subgraph is counted in M. \\Consequently, F$_{21}$= 12$\displaystyle{\sum_{j\neq i}}a_{ij}^{(2)}(d_{j}-a_{ij}-1)(a_{ij}^{(2)}a_{ij})- $ $\frac{3}{4}$ F$_{7}-$  F$_{10}$.\\
\unitlength=0.9mm
\special{em:linewidth 0.4pt}
\linethickness{0.4pt}
\begin{picture}(83.33,35)
\put(55.67,2){\line(-1,2){7}}
\put(73.33,2){\line(1,2){7}}
\put(80.5,16.5){\line(-1,1){16}}
\put(48.90,16.5){\line(1,1){16}}
\put(55.67,2){\line(1,0){18}}
\put(48.90,16){\line(1,0){31}}
\put(48.90,16){\circle*{1.33}}
\put(64.7,32){\circle*{1.33}}
\put(80.5,16.5){\circle*{1.33}}
\put(55.67,2){\circle*{1.33}}
\put(73.33,2){\circle*{1.33}}
\put(93.33,-8){\makebox(0,0)[cc]{Fig 50}}
\put(92.67,-2){\makebox(0,0)[cc]{(b)}}
\put(64.67,-2){\makebox(0,0)[cc]{($a$)}}
\put(103.33,2){\line(1,2){10}}
\put(82.67,2){\line(1,2){10}}
\put(103.33,2){\line(-1,2){10}}
\put(92.90,22){\line(1,0){20.67}}
\put(82.67,2){\line(1,0){20.67}}
\put(92.90,22){\circle*{1.33}}
\put(82.67,2){\circle*{1.33}}
\put(103.33,2){\circle*{1.33}}
\put(113.33,22){\circle*{1.33}}

\put(122.67,-2){\makebox(0,0)[cc]{(c)}}
\put(133.33,2){\line(1,2){10}}
\put(112.67,2){\line(1,2){10}}
\put(133.33,2){\line(-1,2){10}}
\put(122.90,22){\line(1,0){20.67}}
\put(112.67,2){\line(1,0){20.67}}
\put(122.90,22){\circle*{1.33}}
\put(112.67,2){\circle*{1.33}}
\put(133.33,2){\circle*{1.33}}
\put(143.33,22){\circle*{1.33}}
\put(46.47,17){\makebox(0,0)[cc]{$v_{i}$}}
\put(92.47,25){\makebox(0,0)[cc]{$v_{i}$}}
\put(142.47,25){\makebox(0,0)[cc]{$v_{i}$}}
\end{picture}\\\\\\\\
\noindent {\bf Case 22:}
For the configuration of Fig 51(a), N= $6$, M=$\frac{1}{2}$ $\displaystyle{\sum_{j\neq i}}$$a_{ij}^{(2)}(d_{j}-a_{ij}-1)(a_{jj}^{(3)}a_{ij})$ (See theorem 1.7). Let P$_{1}$ denote the number of all subgraphs of G that have the same configuration as the graph of Figure 51(b) and are counted in M. Thus P$_{1}$= $2\times(\frac{1}{24}$ F$_{10})$, where $\frac{1}{24}$ F$_{10}$ is the number of subgraphs of G that have the same configuration as the graph of Figure 51(b) and $2$ is the number of times that this subgraph is counted in M. Let P$_{2}$ denote the number of all subgraphs of G that have the same configuration as the graph of Figure 51(c) and are counted in M. Thus P$_{2}$= $6\times(\frac{1}{12}$ F$_{9})$, where $\frac{1}{12}$ F$_{9}$ is the number of subgraphs of G that have the same configuration as the graph of Figure 51(c) and $6$ is the number of times that this subgraph is counted in M. Let P$_{3}$ denotes the number of all subgraphs of G that have the same configuration as the graph of Figure 51(d) and are counted in M. Thus P$_{3}$= $1\times(\frac{1}{12}$ F$_{21})$, where $\frac{1}{12}$ F$_{21}$ is the number of subgraphs of G that have the same configuration as the graph of Figure 51(d) and $1$ is the number of times that this subgraph is counted in M. Let P$_{4}$ denote the number of all subgraphs of G that have the same configuration as the graph of Figure 51(e) and are counted in M. Thus P$_{4}$= $2\times(\frac{1}{32}$ F$_{7})$, where $\frac{1}{32}$ F$_{7}$ is the number of subgraphs of G that have the same configuration as the graph of Figure 51(e) and $2$ is the number of times that this subgraph is counted in M. Let P$_{5}$ denote the number of all subgraphs of G that have the same configuration as the graph of Figure 51(f) and are counted in M. Thus P$_{5}$= $1\times(\frac{1}{4}$ F$_{29})$, where $\frac{1}{4}$ F$_{29}$ is the number of subgraphs of G that have the same configuration as the graph of Figure 51(f) and $1$ is the number of times that this subgraph is counted in M. Consequently, F$_{22}$= 3 $\displaystyle{\sum_{j\neq i}}$$a_{ij}^{(2)}(d_{j}-a_{ij}-1)(a_{jj}^{(3)}a_{ij})- $ $\frac{1}{2}$ F$_{10}-$ $3$ F$_{9}-$ $\frac{1}{2}$ F$_{21}-$ $\frac{3}{8}$ F$_{7}-$ $\frac{3}{2}$ F$_{29}$.\\
\unitlength=0.9mm
\special{em:linewidth 0.4pt}
\linethickness{0.4pt}
\begin{picture}(83.33,35)
\put(0.67,2){\line(-1,2){7}}
\put(18.33,2){\line(1,2){7}}
\put(25.5,16.5){\line(-1,1){16}}
\put(-5.90,16.5){\line(1,1){16}}
\put(0.67,2){\line(1,0){18}}
\put(-5.90,16){\line(1,0){31}}
\put(-5.90,16){\circle*{1.33}}
\put(9.7,32){\circle*{1.33}}
\put(25.5,16.5){\circle*{1.33}}
\put(0.67,2){\circle*{1.33}}
\put(18.33,2){\circle*{1.33}}
\put(95.33,-6){\makebox(0,0)[cc]{Fig 51}}
\put(37.67,-2){\makebox(0,0)[cc]{(b)}}
\put(9.67,-2){\makebox(0,0)[cc]{($a$)}}
\put(48.33,2){\line(1,2){10}}
\put(27.67,2){\line(1,2){10}}
\put(48.33,2){\line(-1,2){10}}
\put(37.90,22){\line(1,0){20.67}}
\put(27.67,2){\line(1,0){20.67}}
\put(37.90,22){\circle*{1.33}}
\put(27.67,2){\circle*{1.33}}
\put(48.33,2){\circle*{1.33}}
\put(58.33,22){\circle*{1.33}}

\put(67.1,2){\line(1,0){20}}
\put(67.1,2){\line(0,1){20}}
\put(87.1,2){\line(0,1){20}}
\put(87.1,2){\line(1,3){10}}
\put(67.1,22){\line(1,0){20}}
\put(67.1,22){\line(3,1){30}}
\put(67.1,2){\circle*{1.33}}
\put(87.1,2){\circle*{1.33}}
\put(87.1,22){\circle*{1.33}}
\put(67.1,22){\circle*{1.33}}
\put(96.73,31.6){\circle*{1.33}}
\put(87.1,2){\line(-1,1){20}}
\put(25,2){\makebox(0,0)[cc]{$v_{i}$}}
\put(-2,2){\makebox(0,0)[cc]{$v_{i}$}}
\put(64.5,2){\makebox(0,0)[cc]{$v_{i}$}}
\put(77.67,-2){\makebox(0,0)[cc]{(c)}}

\put(105.67,2){\line(-1,2){7}}
\put(123.33,2){\line(1,2){7}}
\put(130.5,16.5){\line(-1,1){16}}
\put(98.90,16.5){\line(1,1){16}}
\put(105.67,2){\line(1,0){18}}
\put(98.90,16){\line(1,0){31}}
\put(98.90,16){\circle*{1.33}}
\put(114.7,32){\circle*{1.33}}
\put(130.5,16.5){\circle*{1.33}}
\put(105.67,2){\circle*{1.33}}
\put(123.33,2){\circle*{1.33}}
\put(95.97,16){\makebox(0,0)[cc]{$v_{i}$}}
\put(141.97,24){\makebox(0,0)[cc]{$v_{i}$}}
\put(142.67,-2){\makebox(0,0)[cc]{(e)}}
\put(114.67,-2){\makebox(0,0)[cc]{(d)}}
\put(153.33,2){\line(1,2){10}}
\put(132.67,2){\line(1,2){10}}
\put(153.33,2){\line(-1,2){10}}
\put(142.90,22){\line(1,0){20.67}}
\put(132.67,2){\line(1,0){20.67}}
\put(142.90,22){\circle*{1.33}}
\put(132.67,2){\circle*{1.33}}
\put(153.33,2){\circle*{1.33}}
\put(163.33,22){\circle*{1.33}}

\put(172.67,2){\line(1,0){20}}
\put(172.67,2){\line(0,1){20}}
\put(192.67,2){\line(0,1){20}}
\put(172.67,22){\line(1,0){20}}
\put(192.67,22){\line(0,1){10}}
\put(192.67,22){\line(1,0){10}}
\put(192.67,32){\line(1,-1){10}}
\put(172.67,2){\circle*{1.33}}
\put(192.67,2){\circle*{1.33}}
\put(202.67,22){\circle*{1.33}}
\put(192.67,32){\circle*{1.33}}
\put(192.67,22){\circle*{1.33}}
\put(172.67,22){\circle*{1.33}}
\put(172,24){\makebox(0,0)[cc]{$v_{i}$}}
\put(182.90,-4){\makebox(0,0)[cc]{$(f)$}}
\end{picture}\\\\\\\\
\noindent {\bf Case 23:}
For the configuration of Figure 52(a), N= 2, M= $\frac{1}{2}$ $ [\displaystyle{\sum_{j\neq i}}[a_{jj}^{(5)} - 5 a_{jj}^{(3)} - 2 ( d_{j} - 2 ) a_{jj}^{(3)}-2$ $\displaystyle{\sum_{k\neq j}}$ $a_{jk}^{(2)} a_{jk}(d_{k}-2) - 2$ $\displaystyle{\sum_{k\neq j}}$ $a_{jk}(\frac{1}{2} a_{kk}^{(3)}-a_{jk}a_{jk}^{(2)})]a_{ij}]$ (See Theorem 1.11). Let P$_{1}$  denote the number of all subgraphs of G that have the same configuration as the graph of Fig 52(b) and are counted in M. Thus P$_{1}=$ $2\times(\frac{1}{14}$ F$_{19})$, where $\frac{1}{14}$ F$_{19}$ is the number of subgraphs of G that have the same configuration as the graph of Figure 52(b) and $2$ is the number of times that this subgraph is counted in M. Let P$_{2}$ denote the number of all subgraphs of G that have the same configuration as the graph of Figure 52(c) and are counted in M. Thus P$_{2}$= $1\times(\frac{1}{12}$ F$_{21})$, where $\frac{1}{12}$ F$_{21}$ is the number of subgraphs of G that have the same configuration as the graph of Figure 52(c) and $1$ is the number of times that this subgraph is counted in M. Consequently, F$_{23}$= $\displaystyle{\sum_{j\neq i}}[a_{jj}^{(5)} - 5 a_{jj}^{(3)} - 2 ( d_{j} - 2 ) a_{jj}^{(3)}-2$ $\displaystyle{\sum_{k\neq j}}$ $a_{jk}^{(2)} a_{jk}(d_{k}-2) - 2$ $\displaystyle{\sum_{k\neq j}}$ $a_{jk}(\frac{1}{2} a_{kk}^{(3)}-a_{jk}a_{jk}^{(2)})]a_{ij}-  \frac{2}{7}$ F$_{19}- \frac{1}{6}$ F$_{21}$.\\\\
\unitlength=0.9mm
\special{em:linewidth 0.4pt}
\linethickness{0.4pt}
\begin{picture}(83.33,40)
\put(45.67,2){\line(-1,2){7}}
\put(63.33,2){\line(1,2){7}}
\put(70.5,16.5){\line(-1,1){16}}
\put(38.90,16.5){\line(1,1){16}}
\put(45.67,2){\line(1,0){18}}
\put(54.7,32){\line(1,0){20}}
\put(74.7,32){\circle*{1.33}}
\put(38.90,16){\circle*{1.33}}
\put(54.7,32){\circle*{1.33}}
\put(70.5,16.5){\circle*{1.33}}
\put(45.67,2){\circle*{1.33}}
\put(63.33,2){\circle*{1.33}}
\put(85.67,2){\line(-1,2){7}}
\put(103.33,2){\line(1,2){7}}
\put(110.5,16.5){\line(-1,1){16}}
\put(78.90,16.5){\line(1,1){16}}
\put(85.67,2){\line(1,0){18}}
\put(78.90,16){\circle*{1.33}}
\put(94.7,32){\circle*{1.33}}
\put(110.5,16.5){\circle*{1.33}}
\put(85.4,2){\circle*{1.33}}
\put(103.33,2){\circle*{1.33}}

\put(125.67,2){\line(-1,2){7}}
\put(143.33,2){\line(1,2){7}}
\put(150.5,16.5){\line(-1,1){16}}
\put(118.90,16.5){\line(1,1){16}}
\put(125.67,2){\line(1,0){18}}
\put(134.7,32){\line(-1,-3){10}}
\put(118.90,16){\circle*{1.33}}
\put(134.7,32){\circle*{1.33}}
\put(150.5,16.5){\circle*{1.33}}
\put(125.4,2){\circle*{1.33}}
\put(143.33,2){\circle*{1.33}}

\put(94.33,-8){\makebox(0,0)[cc]{Fig 52}}
\put(133.67,-2){\makebox(0,0)[cc]{(c)}}
\put(93.67,-2){\makebox(0,0)[cc]{(b)}}
\put(54.67,-2){\makebox(0,0)[cc]{($a$)}}
\put(75,35){\makebox(0,0)[cc]{$v_{i}$}}
\put(95,35){\makebox(0,0)[cc]{$v_{i}$}}
\put(135,35){\makebox(0,0)[cc]{$v_{i}$}}
\end{picture}\\\\\\\\
\noindent {\bf Case 24:}
For the configuration of Figure 53(a), N= 4, M= $\frac{1}{2}$ $[(a_{ii}^{(5)} - 5 a_{ii}^{(3)} - 2 ( d_{i} - 2 ) a_{ii}^{(3)})(d_{i}-2)- 2\displaystyle{\sum_{j\neq i }}a_{ij}^{(2)} a_{ij}(d_{j}-2)(d_{i}-2) - 2\displaystyle{\sum_{j\neq i}}a_{ij}(d_{i}-2)(\frac{1}{2} a_{jj}^{(3)}-a_{ij}a_{ij}^{(2)})]$ (See Theorem 1.11). Let P$_{1}$  denote the number of all subgraphs of G that have the same configuration as the graph of Fig 53(b) and are counted in M. Thus P$_{1}=1\times(\frac{1}{12}$ F$_{21})$, where $\frac{1}{12}$ F$_{21}$ is the number of subgraphs of G that have the same configuration as the graph of Figure 53(b) and $1$ is the number of times that this Fig is counted in M. Consequently, \\F$_{24}$= 2$[(a_{ii}^{(5)} - 5 a_{ii}^{(3)} - 2 ( d_{i} - 2 ) a_{ii}^{(3)})(d_{i}-2)- 2\displaystyle{\sum_{j\neq i }}a_{ij}^{(2)} a_{ij}(d_{j}-2)(d_{i}-2) - 2\displaystyle{\sum_{j\neq i}}a_{ij}(d_{i}-2)(\frac{1}{2} a_{jj}^{(3)}-a_{ij}a_{ij}^{(2)})]-$ $\frac{1}{3}$ F$_{21}$.\\
\unitlength=0.9mm
\special{em:linewidth 0.4pt}
\linethickness{0.4pt}
\begin{picture}(83.33,36)
\put(65.67,2){\line(-1,2){7}}
\put(83.33,2){\line(1,2){7}}
\put(90.5,16.5){\line(-1,1){16}}
\put(58.90,16.5){\line(1,1){16}}
\put(65.67,2){\line(1,0){18}}
\put(74.7,32){\line(1,0){20}}
\put(94.7,32){\circle*{1.33}}
\put(58.90,16){\circle*{1.33}}
\put(74.7,32){\circle*{1.33}}
\put(90.5,16.5){\circle*{1.33}}
\put(65.67,2){\circle*{1.33}}
\put(83.33,2){\circle*{1.33}}
\put(105.67,2){\line(-1,2){7}}
\put(123.33,2){\line(1,2){7}}
\put(130.5,16.5){\line(-1,1){16}}
\put(98.90,16.5){\line(1,1){16}}
\put(105.67,2){\line(1,0){18}}
\put(114.7,32){\line(-1,-3){10}}
\put(98.90,16){\circle*{1.33}}
\put(114.7,32){\circle*{1.33}}
\put(130.5,16.5){\circle*{1.33}}
\put(105.4,2){\circle*{1.33}}
\put(123.33,2){\circle*{1.33}}
\put(95.33,-6){\makebox(0,0)[cc]{Fig 53}}
\put(112.67,-2){\makebox(0,0)[cc]{(b)}}
\put(74.67,-2){\makebox(0,0)[cc]{($a$)}}
\put(75,35){\makebox(0,0)[cc]{$v_{i}$}}
\put(115,35){\makebox(0,0)[cc]{$v_{i}$}}
\end{picture}\\\\\\
\noindent {\bf Case 25:}
For the configuration of Figure 54(a), N= 2, M= $\displaystyle{\sum_{j\neq i}}[a_{ij}^{(4)}-(d_{i}+d_{j}-3a_{ij})a_{ij}^{(2)}-(a_{ii}^{(3)}+a_{jj}^{(3)}+2\left(\begin{array}{c}d_{j}-1\\2\end{array}\right))a_{ij}](d_{j}-2)a_{ij}$ (See Theorem 1.8).  Let P$_{1}$  denote the number of all subgraphs of G that have the same configuration as the graph of Fig 54(b) and are counted in M. Thus P$_{1}=2\times(\frac{1}{6}$ F$_{20})$, where $\frac{1}{6}$ F$_{20}$ is the number of subgraphs of G that have the same configuration as the graph of Figure 54(b) and $2$ is the number of times that this subgraph is counted in M. Let P$_{2}$ denote the number all subgraphs of G that have the same configuration as the graph of Figure 54(c) and are counted in M. Thus P$_{2}$= $1\times(\frac{1}{6}$ F$_{22})$, where $\frac{1}{6}$ F$_{22}$ is the number of subgraphs of G that have the same configuration as the graph of Figure 54(c) and $1$ is the number of times that this subgraph is counted in M. Consequently, F$_{25}= $ 2 $\displaystyle{\sum_{j\neq i}}[a_{ij}^{(4)}-(d_{i}+d_{j}-3a_{ij})a_{ij}^{(2)}-(a_{ii}^{(3)}+a_{jj}^{(3)}+2\left(\begin{array}{c}d_{j}-1\\2\end{array}\right))a_{ij}](d_{j}-2)a_{ij}- $ $\frac{2}{3}$ F$_{20}-$ $\frac{1}{3}$ F$_{22}$.\\
\unitlength=0.9mm
\special{em:linewidth 0.4pt}
\linethickness{0.4pt}
\begin{picture}(83.33,35)
\put(45.67,2){\line(-1,2){7}}
\put(63.33,2){\line(1,2){7}}
\put(70.5,16.5){\line(-1,1){16}}
\put(38.90,16.5){\line(1,1){16}}
\put(45.67,2){\line(1,0){18}}
\put(54.7,32){\line(1,0){20}}
\put(74.7,32){\circle*{1.33}}
\put(38.90,16){\circle*{1.33}}
\put(54.7,32){\circle*{1.33}}
\put(70.5,16.5){\circle*{1.33}}
\put(45.67,2){\circle*{1.33}}
\put(63.33,2){\circle*{1.33}}
\put(85.67,2){\line(-1,2){7}}
\put(103.33,2){\line(1,2){7}}
\put(110.5,16.5){\line(-1,1){16}}
\put(78.90,16.5){\line(1,1){16}}
\put(85.67,2){\line(1,0){18}}
\put(94.7,32){\line(-1,-3){10}}
\put(78.90,16){\circle*{1.33}}
\put(94.7,32){\circle*{1.33}}
\put(110.5,16.5){\circle*{1.33}}
\put(85.4,2){\circle*{1.33}}
\put(103.33,2){\circle*{1.33}}

\put(125.67,2){\line(-1,2){7}}
\put(143.33,2){\line(1,2){7}}
\put(150.5,16.5){\line(-1,1){16}}
\put(118.90,16.5){\line(1,1){16}}
\put(125.67,2){\line(1,0){18}}
\put(135.33,32){\line(-1,-3){10}}
\put(118.90,16){\circle*{1.33}}
\put(135,32){\circle*{1.33}}
\put(150.5,16.5){\circle*{1.33}}
\put(125.4,2){\circle*{1.33}}
\put(143.33,2){\circle*{1.33}}

\put(95.33,-8){\makebox(0,0)[cc]{Fig 54}}
\put(133.67,-2){\makebox(0,0)[cc]{(c)}}
\put(93.67,-2){\makebox(0,0)[cc]{(b)}}
\put(54.67,-2){\makebox(0,0)[cc]{($a$)}}
\put(36,16.5){\makebox(0,0)[cc]{$v_{i}$}}
\put(76,16.5){\makebox(0,0)[cc]{$v_{i}$}}
\put(153.4,16.5){\makebox(0,0)[cc]{$v_{i}$}}
\end{picture}\\\\\\\\
\noindent {\bf Case 26:}
For the configuration of Figure 55(a), N= 2, M= $ \displaystyle{\sum_{j\neq i}}[(a_{ij}^{(3)} - (d_{i}+ d_{j}- 1)a_{ij})a_{ij}^{(2)}- \displaystyle{\sum_{k\neq i,j}}(a_{ik}^{(2)}a_{ik}- a_{ik}a_{ij}a_{jk})a_{jk}- \displaystyle{\sum_{k\neq i,j}}(a_{jk}^{(2)}a_{jk}- a_{ik}a_{ij}a_{jk})a_{ik}](d_{j}- 2)$. Let P$_{1}$  denote the number of all subgraphs of G that have the same configuration as the graph of Fig 55(b) and are counted in M. Thus P$_{1}=1\times(\frac{1}{12}$ F$_{21})$, where $\frac{1}{12}$ F$_{21}$ is the number of subgraphs of G that have the same configuration as the graph of Fig 55(b) and 1 is the number of times that this subgraph is counted in M. Let P$_{2}$  denote the number of all subgraphs of G that have the same configuration as the graph of Fig 55(c) and are counted in M. Thus P$_{2}=1\times(\frac{1}{6}$ F$_{22})$, where $\frac{1}{6}$ F$_{22}$ is the number of subgraphs of G that have the same configuration as the graph of Fig 55(c) and 1 is the number of times that this subgraph is counted in M. Consequently, \\F= 2$ \displaystyle{\sum_{j\neq i}}[(a_{ij}^{(3)} - (d_{i}+ d_{j}- 1)a_{ij})a_{ij}^{(2)}- \displaystyle{\sum_{k\neq i,j}}(a_{ik}^{(2)}a_{ik}- a_{ik}a_{ij}a_{jk})a_{jk}- \displaystyle{\sum_{k\neq i,j}}(a_{jk}^{(2)}a_{jk}- a_{ik}a_{ij}a_{jk})a_{ik}](d_{j}- 2)-$ $\frac{1}{6}$ F$_{21}\\- \frac{1}{3}$ F$_{22}$ .\\\\\\
\unitlength=0.9mm
\special{em:linewidth 0.4pt}
\linethickness{0.4pt}
\begin{picture}(83.33,27)
\put(55.67,2){\line(-1,2){7}}
\put(73.33,2){\line(1,2){7}}
\put(80.5,16.5){\line(-1,1){16}}
\put(48.90,16.5){\line(1,1){16}}
\put(55.67,2){\line(1,0){18}}
\put(64.7,32){\line(1,0){20}}
\put(84.7,32){\circle*{1.33}}
\put(48.90,16){\circle*{1.33}}
\put(64.7,32){\circle*{1.33}}
\put(80.5,16.5){\circle*{1.33}}
\put(55.67,2){\circle*{1.33}}
\put(73.33,2){\circle*{1.33}}

\put(135.67,2){\line(-1,2){7}}
\put(153.33,2){\line(1,2){7}}
\put(160.5,16.5){\line(-1,1){16}}
\put(128.90,16.5){\line(1,1){16}}
\put(135.67,2){\line(1,0){18}}
\put(128.90,16){\line(1,0){31}}
\put(128.90,16){\circle*{1.33}}
\put(144.7,32){\circle*{1.33}}
\put(160.5,16.5){\circle*{1.33}}
\put(135.4,2){\circle*{1.33}}
\put(153.33,2){\circle*{1.33}}
\put(103.33,-8){\makebox(0,0)[cc]{Fig 55}}
\put(143.67,-2){\makebox(0,0)[cc]{(c)}}
\put(104.67,-2){\makebox(0,0)[cc]{(b)}}
\put(63.67,-2){\makebox(0,0)[cc]{($a$)}}

\put(95.67,2){\line(-1,2){7}}
\put(113.33,2){\line(1,2){7}}
\put(120.5,16.5){\line(-1,1){16}}
\put(88.90,16.5){\line(1,1){16}}
\put(95.67,2){\line(1,0){18}}
\put(88.90,16){\line(1,0){31}}
\put(88.90,16){\circle*{1.33}}
\put(104.7,32){\circle*{1.33}}
\put(120.5,16.5){\circle*{1.33}}
\put(95.67,2){\circle*{1.33}}
\put(113.33,2){\circle*{1.33}}
\put(53,2){\makebox(0,0)[cc]{$v_{i}$}}
\put(86,16){\makebox(0,0)[cc]{$v_{i}$}}
\put(133,2){\makebox(0,0)[cc]{$v_{i}$}}
\end{picture}\\\\\\\\
\noindent {\bf Case 27:}
For the configuration of Figure $56( a )$, N= $4$, M= $\displaystyle{\sum_{j\neq i}}[a_{ij}^{(2)}a_{ij} \times \displaystyle{\sum_{k\neq i,j}}  \left(\begin{array}{c}a_{ij}^{(2)}\\2\end{array}\right)]$. Let Let P$_{1}$  denote the number of all subgraphs of G that have the same configuration as the graph of Fig 56(b) and are counted in M. Thus P$_{1}=1\times(\frac{1}{12}$ F$_{21})$, where $\frac{1}{12}$ F$_{21}$ is the number of subgraphs of G that have the same configuration as the graph of Fig 56(b) and 1 is the number of times that this subgraph is counted in M. Let P$_{2}$  denote the number of all subgraphs of G that have the same configuration as the graph of Fig 56(c) and are counted in M. Thus P$_{2}=2\times(\frac{1}{24}$ F$_{10})$, where $\frac{1}{24}$ F$_{10}$ is the number of subgraphs of G that have the same configuration as the graph of Fig 56(c) and 2 is the number of times that this subgraph is counted in M. Let P$_{3}$  denote the number of all subgraphs of G that have the same configuration as the graph of Fig 56(d) and are counted in M. Thus P$_{3}=2\times(\frac{1}{12}$ F$_{9})$, where $\frac{1}{12}$ F$_{9}$ is the number of subgraphs of G that have the same configuration as the graph of Fig 56(d) and 2 is the number of times that this subgraph is counted in M. Let P$_{4}$  denote the number of all subgraphs of G that have the same configuration as the graph of Fig 56(e) and are counted in M. Thus P$_{4}=2\times(\frac{1}{6}$ F$_{20})$, where $\frac{1}{6}$ F$_{20}$ is the number of subgraphs of G that have the same configuration as the graph of Fig 56(e) and 2 is the number of times that this subgraph is counted in M. Let P$_{5}$  denote the number of all subgraphs of G that have the same configuration as the graph of Fig 56(f) and are counted in M. Thus P$_{5}=2\times(\frac{1}{32}$ F$_{7})$, where $\frac{1}{32}$ F$_{7}$ is the number of subgraphs of G that have the same configuration as the graph of Fig 56(f) and 2 is the number of times that this subgraph is counted in M. Consequently, F= 4$\displaystyle{\sum_{j\neq i}}[a_{ij}^{(2)}a_{ij} \times \displaystyle{\sum_{k\neq i,j}}  \left(\begin{array}{c}a_{ij}^{(2)}\\2\end{array}\right)]- $ $\frac{1}{3}$ F$_{21}-$ $\frac{1}{3}$ F$_{10}-$ $\frac{2}{3}$ F$_{9}-$ $\frac{4}{3}$ F$_{20}-$ $\frac{1}{4}$ F$_{7}$.\\
\unitlength=0.9mm
\special{em:linewidth 0.4pt}
\linethickness{0.4pt}
\begin{picture}(83.33,35)
\put(42.67,2){\line(1,0){20}}
\put(42.67,2){\line(0,1){20}}
\put(62.67,2){\line(0,1){20}}
\put(62.67,22){\line(-1,1){10}}
\put(42.67,22){\line(1,1){10}}
\put(42.67,22){\line(1,0){20}}
\put(42.67,2){\circle*{1.33}}
\put(52.67,32){\circle*{1.33}}
\put(62.67,2){\circle*{1.33}}
\put(62.67,22){\circle*{1.33}}
\put(42.67,22){\circle*{1.33}}
\put(40,22){\makebox(0,0)[cc]{$v_{i}$}}
\put(51.90,-4){\makebox(0,0)[cc]{$(b)$}}
\put(141.90,-4){\makebox(0,0)[cc]{$(e)$}}

\put(12.67,2){\line(1,0){20}}
\put(12.67,2){\line(0,1){20}}
\put(32.67,2){\line(0,1){20}}
\put(12.67,22){\line(1,0){20}}
\put(12.67,22){\line(0,1){10}}
\put(12.67,22){\line(-1,0){10}}
\put(12.67,32){\line(-1,-1){10}}
\put(32.67,2){\circle*{1.33}}
\put(32.67,2){\circle*{1.33}}
\put(12.67,2){\circle*{1.33}}
\put(2.67,22){\circle*{1.33}}
\put(12.67,32){\circle*{1.33}}
\put(32.67,22){\circle*{1.33}}
\put(12.67,22){\circle*{1.33}}
\put(16,32){\makebox(0,0)[cc]{$v_{i}$}}
\put(22.90,-4){\makebox(0,0)[cc]{$($a$)$}}

\put(72.67,2){\line(1,0){20}}
\put(92.67,2){\line(-1,1){20}}
\put(72.67,2){\line(0,1){20}}
\put(92.67,2){\line(0,1){20}}
\put(72.67,22){\line(1,0){20}}
\put(72.67,2){\circle*{1.33}}
\put(92.67,2){\circle*{1.33}}
\put(92.67,22){\circle*{1.33}}
\put(72.67,22){\circle*{1.33}}
\put(93,24){\makebox(0,0)[cc]{$v_{i}$}}
\put(81.90,-4){\makebox(0,0)[cc]{$(c)$}}

\put(102.67,2){\line(1,0){20}}
\put(122.67,2){\line(-1,1){20}}
\put(102.67,2){\line(0,1){20}}
\put(122.67,2){\line(0,1){20}}
\put(122.67,2){\line(1,3){10}}
\put(102.67,22){\line(1,0){20}}
\put(102.67,22){\line(3,1){30}}
\put(102.67,2){\circle*{1.33}}
\put(122.67,2){\circle*{1.33}}
\put(122.67,22){\circle*{1.33}}
\put(102.67,22){\circle*{1.33}}
\put(132.3,31.6){\circle*{1.33}}
\put(132,34){\makebox(0,0)[cc]{$v_{i}$}}
\put(111.90,-4){\makebox(0,0)[cc]{$(d)$}}
\put(96.90,-8){\makebox(0,0)[cc]{Fig 56}}

\put(132.67,2){\line(1,0){20}}
\put(132.67,2){\line(0,1){20}}
\put(152.67,2){\line(0,1){20}}
\put(152.67,22){\line(-1,1){10}}
\put(132.67,22){\line(1,1){10}}
\put(132.67,22){\line(1,0){20}}
\put(132.67,2){\circle*{1.33}}
\put(142.67,32){\circle*{1.33}}
\put(152.67,2){\circle*{1.33}}
\put(152.67,22){\circle*{1.33}}
\put(132.67,22){\circle*{1.33}}
\put(143,34){\makebox(0,0)[cc]{$v_{i}$}}

\put(162.67,2){\line(1,0){20}}
\put(182.67,2){\line(-1,1){20}}
\put(162.67,2){\line(0,1){20}}
\put(182.67,2){\line(0,1){20}}
\put(162.67,22){\line(1,0){20}}
\put(162.67,2){\circle*{1.33}}
\put(182.67,2){\circle*{1.33}}
\put(182.67,22){\circle*{1.33}}
\put(162.67,22){\circle*{1.33}}
\put(160,22){\makebox(0,0)[cc]{$v_{i}$}}
\put(171.90,-4){\makebox(0,0)[cc]{$(f)$}}

\end{picture}\\\\\\\\
\noindent {\bf Case 28:}
For the configuration of Figure $57( a )$, N=8, M= $(\frac{1}{2}a_{ii}^{(3)})\displaystyle{\sum_{j\neq i}}  \left(\begin{array}{c}a_{ij}^{(2)}\\2\end{array}\right)$. Let P$_{1}$  denote the number of all subgraphs of G that have the same configuration as the graph of Fig 57(b) and are counted in M. Thus P$_{1}=1\times(\frac{1}{12}$ F$_{21})$, where $\frac{1}{12}$ F$_{21}$ is the number of subgraphs of G that have the same configuration as the graph of Fig 57(b) and 1 is the number of times that this subgraph is counted in M. Let P$_{2}$ denote the number of all subgraphs of G that have the same configuration as the graph of Fig 57(c) and are counted in M. Thus P$_{2}=1\times(\frac{1}{24}$ F$_{10})$, where $\frac{1}{24}$ F$_{10}$ is the number of subgraphs of G that have the same configuration as the graph of Fig 57(c) and 1 is the number of times that this subgraph is counted in M. Let P$_{3}$  denote the number of all subgraphs of G that have the same configuration as the graph of Fig 57(d) and are counted in M. Thus P$_{3}=3\times(\frac{1}{24}$ F$_{8})$, where $\frac{1}{24}$ F$_{8}$ is the number of subgraphs of G that have the same configuration as the graph of Fig 57(d) and 3 is the number of times that this subgraph is counted in M. Let P$_{4}$  denote the number of all subgraphs of G that have the same configuration as the graph of Fig 57(e) and are counted in M. Thus P$_{4}=2\times(\frac{1}{32}$ F$_{7})$, where $\frac{1}{32}$ F$_{7}$ is the number of subgraphs of G that have the same configuration as the graph of Fig 57(e) and 2 is the number of times that this subgraph is counted in M. Consequently, F= $(4 a_{ii}^{(3)})\displaystyle{\sum_{j\neq i}}  \left(\begin{array}{c}a_{ij}^{(2)}\\2\end{array}\right)-$ $\frac{2}{3}$ F$_{21}-$ $\frac{1}{3}$ F$_{10}-$ F$_{8}-$ $\frac{1}{2}$ F$_{7}$.\\
\unitlength=0.9mm
\special{em:linewidth 0.4pt}
\linethickness{0.4pt}
\begin{picture}(83.33,35)
\put(52.67,2){\line(1,0){20}}
\put(52.67,2){\line(0,1){20}}
\put(72.67,2){\line(0,1){20}}
\put(72.67,22){\line(-1,1){10}}
\put(52.67,22){\line(1,1){10}}
\put(52.67,22){\line(1,0){20}}
\put(52.67,2){\circle*{1.33}}
\put(62.67,32){\circle*{1.33}}
\put(72.67,2){\circle*{1.33}}
\put(72.67,22){\circle*{1.33}}
\put(52.67,22){\circle*{1.33}}
\put(50,22){\makebox(0,0)[cc]{$v_{i}$}}
\put(61.90,-4){\makebox(0,0)[cc]{$(b)$}}

\put(22.67,2){\line(1,0){20}}
\put(22.67,2){\line(0,1){20}}
\put(42.67,2){\line(0,1){20}}
\put(22.67,22){\line(1,0){20}}
\put(22.67,22){\line(0,1){10}}
\put(22.67,22){\line(-1,0){10}}
\put(22.67,32){\line(-1,-1){10}}
\put(42.67,2){\circle*{1.33}}
\put(42.67,2){\circle*{1.33}}
\put(22.67,2){\circle*{1.33}}
\put(12.67,22){\circle*{1.33}}
\put(22.67,32){\circle*{1.33}}
\put(42.67,22){\circle*{1.33}}
\put(22.67,22){\circle*{1.33}}
\put(26,24){\makebox(0,0)[cc]{$v_{i}$}}
\put(32.90,-4){\makebox(0,0)[cc]{$($a$)$}}

\put(82.67,2){\line(1,0){20}}
\put(102.67,2){\line(-1,1){20}}
\put(82.67,2){\line(0,1){20}}
\put(102.67,2){\line(0,1){20}}
\put(82.67,22){\line(1,0){20}}
\put(82.67,2){\circle*{1.33}}
\put(102.67,2){\circle*{1.33}}
\put(102.67,22){\circle*{1.33}}
\put(82.67,22){\circle*{1.33}}
\put(103,24){\makebox(0,0)[cc]{$v_{i}$}}
\put(92.90,-4){\makebox(0,0)[cc]{$(c)$}}

\put(112.67,2){\line(1,0){20}}
\put(132.67,2){\line(-1,1){20}}
\put(112.67,2){\line(0,1){20}}
\put(132.67,2){\line(0,1){20}}
\put(132.67,2){\line(1,3){10}}
\put(112.67,22){\line(1,0){20}}
\put(112.67,22){\line(3,1){30}}
\put(112.67,2){\circle*{1.33}}
\put(132.67,2){\circle*{1.33}}
\put(132.67,22){\circle*{1.33}}
\put(112.67,22){\circle*{1.33}}
\put(142.3,31.6){\circle*{1.33}}
\put(112,24){\makebox(0,0)[cc]{$v_{i}$}}
\put(121.90,-4){\makebox(0,0)[cc]{$(d)$}}
\put(93.90,-10){\makebox(0,0)[cc]{Fig 57}}

\put(152.67,2){\line(1,0){20}}
\put(172.67,2){\line(-1,1){20}}
\put(152.67,2){\line(0,1){20}}
\put(172.67,2){\line(0,1){20}}
\put(152.67,22){\line(1,0){20}}
\put(152.67,2){\circle*{1.33}}
\put(172.67,2){\circle*{1.33}}
\put(172.67,22){\circle*{1.33}}
\put(152.67,22){\circle*{1.33}}
\put(150,22){\makebox(0,0)[cc]{$v_{i}$}}
\put(161.90,-4){\makebox(0,0)[cc]{$(e)$}}

\end{picture}\\\\\\\\\\
\noindent {\bf Case 29:}
For the configuration of Figure $58( a )$, N= 4, M= $\frac{1}{2}\displaystyle{\sum_{j\neq i}}[(d_{j}- a_{ij}- 1) a_{jj}^{(3)}a_{ij}^{(2)}a_{ij}]$. Let P$_{1}$  denote the number of all subgraphs of G that have the same configuration as the graph of Fig 58(b) and are counted in M. Thus P$_{1}=1\times(\frac{1}{12}$ F$_{21})$, where $\frac{1}{12}$ F$_{21}$ is the number of subgraphs of G that have the same configuration as the graph of Fig 58(b)  and 1 is the number of times that this subgraph is counted in M. Let P$_{2}$  denote the number of all subgraphs of G that have the same configuration as the graph of Fig 58(c) and are counted in M. Thus P$_{2}=4\times(\frac{1}{24}$ F$_{10})$, where $\frac{1}{24}$ F$_{10}$ is the number of subgraphs of G that have the same configuration as the graph of Fig 58(c) and 4 is the number of times that this subgraph is counted in M. Let P$_{3}$  denote the number of all subgraphs of G that have the same configuration as the graph of Fig 58(d) and are counted in M. Thus P$_{3}=4\times(\frac{1}{12}$ F$_{9})$, where $\frac{1}{12}$ F$_{9}$ is the number of subgraphs of G that have the same configuration as the graph of Fig 58(d) and 4 is the number of times that this subgraph is counted in M. Let P$_{4}$  denote the number of all subgraphs of G that have the same configuration as the graph of Fig 58(e) and are counted in M. Thus P$_{4}=1\times(\frac{1}{6}$ F$_{22})$, where $\frac{1}{6}$ F$_{22}$ is the number of subgraphs of G that have the same configuration as the graph of Fig 58(e) and 1 is the number of times that this subgraph is counted in M. Let P$_{5}$  denote the number of all subgraphs of G that have the same configuration as the graph of Fig 58(f) and are counted in M. Thus P$_{5}=2\times(\frac{1}{32}$ F$_{7})$, where $\frac{1}{32}$ F$_{7}$ is the number of subgraphs of G that have the same configuration as the graph of Fig 58(f) and 2 is the number of times that this subgraph is counted in M. Consequently, F= $2 \displaystyle{\sum_{j\neq i}}[(d_{j}- a_{ij}- 1) a_{jj}^{(3)}a_{ij}^{(2)}a_{ij}]-$$\frac{1}{3}$ F$_{21}-$ $\frac{2}{3}$ F$_{10}-$ $\frac{4}{3}$ F$_{9}-$ $\frac{2}{3}$ F$_{22}-$ $\frac{1}{4}$ F$_{7}$.

\unitlength=0.9mm
\special{em:linewidth 0.4pt}
\linethickness{0.4pt}
\begin{picture}(83.33,35)
\put(42.67,2){\line(1,0){20}}
\put(42.67,2){\line(0,1){20}}
\put(62.67,2){\line(0,1){20}}
\put(62.67,22){\line(-1,1){10}}
\put(42.67,22){\line(1,1){10}}
\put(42.67,22){\line(1,0){20}}
\put(42.67,2){\circle*{1.33}}
\put(52.67,32){\circle*{1.33}}
\put(62.67,2){\circle*{1.33}}
\put(62.67,22){\circle*{1.33}}
\put(42.67,22){\circle*{1.33}}
\put(40,22){\makebox(0,0)[cc]{$v_{i}$}}
\put(51.90,-4){\makebox(0,0)[cc]{$(b)$}}
\put(141.90,-4){\makebox(0,0)[cc]{$(e)$}}

\put(12.67,2){\line(1,0){20}}
\put(12.67,2){\line(0,1){20}}
\put(32.67,2){\line(0,1){20}}
\put(12.67,22){\line(1,0){20}}
\put(12.67,22){\line(0,1){10}}
\put(12.67,22){\line(-1,0){10}}
\put(12.67,32){\line(-1,-1){10}}
\put(32.67,2){\circle*{1.33}}
\put(32.67,2){\circle*{1.33}}
\put(12.67,2){\circle*{1.33}}
\put(2.67,22){\circle*{1.33}}
\put(12.67,32){\circle*{1.33}}
\put(32.67,22){\circle*{1.33}}
\put(12.67,22){\circle*{1.33}}
\put(32,24){\makebox(0,0)[cc]{$v_{i}$}}
\put(22.90,-4){\makebox(0,0)[cc]{$($a$)$}}

\put(72.67,2){\line(1,0){20}}
\put(92.67,2){\line(-1,1){20}}
\put(72.67,2){\line(0,1){20}}
\put(92.67,2){\line(0,1){20}}
\put(72.67,22){\line(1,0){20}}
\put(72.67,2){\circle*{1.33}}
\put(92.67,2){\circle*{1.33}}
\put(92.67,22){\circle*{1.33}}
\put(72.67,22){\circle*{1.33}}
\put(93,24){\makebox(0,0)[cc]{$v_{i}$}}
\put(81.90,-4){\makebox(0,0)[cc]{$(c)$}}

\put(102.67,2){\line(1,0){20}}
\put(122.67,2){\line(-1,1){20}}
\put(102.67,2){\line(0,1){20}}
\put(122.67,2){\line(0,1){20}}
\put(122.67,2){\line(1,3){10}}
\put(102.67,22){\line(1,0){20}}
\put(102.67,22){\line(3,1){30}}
\put(102.67,2){\circle*{1.33}}
\put(122.67,2){\circle*{1.33}}
\put(122.67,22){\circle*{1.33}}
\put(102.67,22){\circle*{1.33}}
\put(132.3,31.6){\circle*{1.33}}
\put(100,2){\makebox(0,0)[cc]{$v_{i}$}}
\put(111.90,-4){\makebox(0,0)[cc]{$(d)$}}
\put(96.90,-8){\makebox(0,0)[cc]{Fig 58}}

\put(132.67,2){\line(1,0){20}}
\put(132.67,2){\line(0,1){20}}
\put(152.67,2){\line(0,1){20}}
\put(152.67,22){\line(-1,1){10}}
\put(132.67,22){\line(1,1){10}}
\put(132.67,22){\line(1,0){20}}
\put(132.67,2){\circle*{1.33}}
\put(142.67,32){\circle*{1.33}}
\put(152.67,2){\circle*{1.33}}
\put(152.67,22){\circle*{1.33}}
\put(132.67,22){\circle*{1.33}}
\put(130,2){\makebox(0,0)[cc]{$v_{i}$}}

\put(162.67,2){\line(1,0){20}}
\put(182.67,2){\line(-1,1){20}}
\put(162.67,2){\line(0,1){20}}
\put(182.67,2){\line(0,1){20}}
\put(162.67,22){\line(1,0){20}}
\put(162.67,2){\circle*{1.33}}
\put(182.67,2){\circle*{1.33}}
\put(182.67,22){\circle*{1.33}}
\put(162.67,22){\circle*{1.33}}
\put(160,22){\makebox(0,0)[cc]{$v_{i}$}}
\put(171.90,-4){\makebox(0,0)[cc]{$(f)$}}

\end{picture}\\\\\\\\
\noindent {\bf Case 30:}
For the configuration of Figure $59( a )$, N= $4$, M= $\frac{1}{2}\displaystyle{\sum_{j\neq i}}\left(\begin{array}{c}a_{ij}^{(2)}\\2\end{array}\right)a_{jj}^{(3)}$. Let P$_{1}$  denote the number of all subgraphs of G that have the same configuration as the graph of Fig 59(b) and are counted in M. Thus P$_{1}=1\times(\frac{1}{6}$ F$_{22})$, where $\frac{1}{6}$ F$_{22}$ is the number of subgraphs of G that have the same configuration as the graph of Fig 59(b) and 1 is the number of times that this subgraph is counted in M. Let P$_{2}$  denote the number of all subgraphs of G that have the same configuration as the graph of Fig 59(c) and are counted in M. Thus P$_{2}=1\times(\frac{1}{24}$ F$_{10})$, where $\frac{1}{24}$ F$_{10}$ is the number of subgraphs of G that have the same configuration as the graph of Fig 59(c) and 1 is the number of times that this subgraph is counted in M. Let P$_{3}$  denote the number of all subgraphs of G that have the same configuration as the graph of Fig 59(d) and are counted in M. Thus P$_{3}=3\times(\frac{1}{24}$ F$_{8})$, where $\frac{1}{24}$ F$_{8}$ is the number of subgraphs of G that have the same configuration as the graph of Fig 59(d) and 3 is the number of times that this subgraph is counted in M. Let P$_{4}$  denote the number of all subgraphs of G that have the same configuration as the graph of Fig 59(e) and are counted in M. Thus P$_{4}=2\times(\frac{1}{32}$ F$_{7})$, where $\frac{1}{32}$ F$_{7}$ is the number of subgraphs of G that have the same configuration as the graph of Fig 59(e) and 2 is the number of times that this subgraph is counted in M. Consequently, F= $2 \displaystyle{\sum_{j\neq i}}\left(\begin{array}{c}a_{ij}^{(2)}\\2\end{array}\right)a_{jj}^{(3)}-$ $\frac{2}{3}$ F$_{22}-$ $\frac{1}{6}$ F$_{10}-$ $\frac{1}{2}$ F$_{8}-$ $\frac{1}{4}$ F$_{7}$.\\
\unitlength=0.9mm
\special{em:linewidth 0.4pt}
\linethickness{0.4pt}
\begin{picture}(83.33,35)
\put(52.67,2){\line(1,0){20}}
\put(52.67,2){\line(0,1){20}}
\put(72.67,2){\line(0,1){20}}
\put(72.67,22){\line(-1,1){10}}
\put(52.67,22){\line(1,1){10}}
\put(52.67,22){\line(1,0){20}}
\put(52.67,2){\circle*{1.33}}
\put(62.67,32){\circle*{1.33}}
\put(72.67,2){\circle*{1.33}}
\put(72.67,22){\circle*{1.33}}
\put(52.67,22){\circle*{1.33}}
\put(75,2){\makebox(0,0)[cc]{$v_{i}$}}
\put(61.90,-4){\makebox(0,0)[cc]{$(b)$}}

\put(22.67,2){\line(1,0){20}}
\put(22.67,2){\line(0,1){20}}
\put(42.67,2){\line(0,1){20}}
\put(22.67,22){\line(1,0){20}}
\put(22.67,22){\line(0,1){10}}
\put(22.67,22){\line(-1,0){10}}
\put(22.67,32){\line(-1,-1){10}}
\put(42.67,2){\circle*{1.33}}
\put(42.67,2){\circle*{1.33}}
\put(22.67,2){\circle*{1.33}}
\put(12.67,22){\circle*{1.33}}
\put(22.67,32){\circle*{1.33}}
\put(42.67,22){\circle*{1.33}}
\put(22.67,22){\circle*{1.33}}
\put(45,2){\makebox(0,0)[cc]{$v_{i}$}}
\put(32.90,-4){\makebox(0,0)[cc]{$($a$)$}}

\put(82.67,2){\line(1,0){20}}
\put(102.67,2){\line(-1,1){20}}
\put(82.67,2){\line(0,1){20}}
\put(102.67,2){\line(0,1){20}}
\put(82.67,22){\line(1,0){20}}
\put(82.67,2){\circle*{1.33}}
\put(102.67,2){\circle*{1.33}}
\put(102.67,22){\circle*{1.33}}
\put(82.67,22){\circle*{1.33}}
\put(103,24){\makebox(0,0)[cc]{$v_{i}$}}
\put(91.90,-4){\makebox(0,0)[cc]{$(c)$}}

\put(112.67,2){\line(1,0){20}}
\put(132.67,2){\line(-1,1){20}}
\put(112.67,2){\line(0,1){20}}
\put(132.67,2){\line(0,1){20}}
\put(132.67,2){\line(1,3){10}}
\put(112.67,22){\line(1,0){20}}
\put(112.67,22){\line(3,1){30}}
\put(112.67,2){\circle*{1.33}}
\put(132.67,2){\circle*{1.33}}
\put(132.67,22){\circle*{1.33}}
\put(112.67,22){\circle*{1.33}}
\put(142.3,31.6){\circle*{1.33}}
\put(112.5,24.5){\makebox(0,0)[cc]{$v_{i}$}}
\put(121.90,-4){\makebox(0,0)[cc]{$(d)$}}
\put(92.90,-10){\makebox(0,0)[cc]{Fig 59}}

\put(152.67,2){\line(1,0){20}}
\put(172.67,2){\line(-1,1){20}}
\put(152.67,2){\line(0,1){20}}
\put(172.67,2){\line(0,1){20}}
\put(152.67,22){\line(1,0){20}}
\put(152.67,2){\circle*{1.33}}
\put(172.67,2){\circle*{1.33}}
\put(172.67,22){\circle*{1.33}}
\put(152.67,22){\circle*{1.33}}
\put(150,22){\makebox(0,0)[cc]{$v_{i}$}}
\put(161.90,-4){\makebox(0,0)[cc]{$(e)$}}

\end{picture}\\\\\\\\
Now we add the values of F$_{n}$ arising from the above cases and determine $x$. Substituting the value of $x$ in \\$\frac{1}{2}(a_{ii}^{(7)}- x)$ and simplifying, we get the number of $7-$cycles each of which contains $a$ specific vertex $v_{i}$ of G.\hfill$\Box$\\
\begin{Example}
In the graph of Fig $29$ we have, $F_{1}= 420,$ $F_{2}= 1020,$ $F_{3}= 1080,$ $F_{4}= 240,$ $F_{5}= 360,$ $F_{6}= 240,$ $F_{7}= 960,$ $F_{8}= 480,$ $F_{9}= 360,$ $F_{10}= 720,$ $F_{11}= 420,$ $F_{12}= 120,$ $F_{13}= 240,$ $F_{14}= 240,$ $F_{15}= 240,$ $F_{16}= 240,$ $F_{17}= 480,$ $F_{18}= 120,$ $F_{19}= 840,$ $F_{20}= 360,$ $F_{21}= 1440,$ $F_{22}= 720,$ $F_{23}= 120,$ $F_{24}= 240,$ $F_{25}= 240,$ $F_{26}= 240,$ $F_{27}= 240,$ $F_{28}= 240,$ $F_{29}= 240,$ $F_{30}= 120,$. So, we have $x=$ $13020$. Consequently, by Theorem $3.3$, the number of $ ~ 7-$cycles each of which contains the vertex $v_{1}$ in the graph of Fig $29$ is $0$.
\end{Example}

\end{document}